\documentclass[12pt, reqno, preprint]{article}

\usepackage[margin= 1.27in]{geometry} %--- turn on this

\usepackage[small]{titlesec} %headings and subheadings
\titleformat{\subsubsection}[runin]{\itshape}{\upshape\thesubsubsection.}{0.5em}{}[.] %JAA

\usepackage{amssymb, amsmath, amsfonts}
\usepackage{hyperref}
\hypersetup{
	colorlinks=true,
	linkcolor=red,
	citecolor=blue,
	filecolor=cyan,      
	urlcolor=magenta}

\usepackage{amsthm}
\newtheoremstyle{remark}{}{}{\upshape}{}{\bfseries}{.}{0.5em}{}
\newtheorem{theorem}{Theorem}[section]

\newtheorem{corollary}[theorem]{Corollary}
\newtheorem{proposition}[theorem]{Proposition}
\newtheorem{definition}[theorem]{Definition}
\theoremstyle{remark}

\newtheorem{remark}[theorem]{Remark}

\usepackage{comment}
\usepackage{bm}

\usepackage[all]{xy}
\usepackage{tikz-cd} 

\newcommand\genfd{{\bm k}}
\newcommand\op{\mathrm{op}}
\newcommand\co{\mathrm{co}}
\newcommand\id{\mathrm{id}}
\newcommand\KK{\mathcal{K}}
\newcommand\ipV{\mathrm{indproVect}}

\newcommand\UgR{U(\mathfrak{g}_R)}
\newcommand\UgL{U(\mathfrak{g}_L)}
\newcommand\Ug{U(\mathfrak{g})}
\newcommand\OO{\mathcal{O}}
\newcommand\Aut{\operatorname{Aut}}
\newcommand\ggf{\mathfrak{g}}

\newcommand\Ad{\operatorname{Ad}}

\newcommand\ad{\operatorname{ad}}
\newcommand\GG{\mathcal{G}}
\newcommand\UU{\mathcal{U}}

\newcommand\End{\operatorname{End}}

\newcommand\gR{\mathfrak{g}_R}
\newcommand\gL{\mathfrak{g}_L}

\newcommand\Diff{\operatorname{Diff}}
\newcommand\GLnk{\mathrm{GL}(n,\genfd)}
\newcommand{\Img}{\operatorname{Im}}

\usepackage{mathabx}
\newcommand\btr{\blacktriangleright}
\newcommand\btl{\blacktriangleleft}
\newcommand\str{\smalltriangleright}

\begin{document}
	
	\title{\Large\bfseries Examples of scalar extension Hopf algebroids over a universal enveloping algebra}
	\author{\normalsize\bfseries M.~Stoji\'c, Z.~\v{Skoda}} 
	\date{\small \today}
%	\address{Department of Mathematics, University of Zagreb, Bijeni\v{c}ka cesta~30, 10000 Zagreb, Croatia}
%	\email{stojic@math.hr} 
	\maketitle

\begin{abstract} 
	The standard pairing between the algebra $\OO(G)$ of regular functions on an affine algebraic group $G$ and the universal enveloping algebra $U(\ggf)$ of the corresponding finite-dimesional Lie algebra $\ggf$ of left invariant derivations of $\OO(G)$ induces a Hopf action of $\OO(G)$ on $U(\ggf)$ which together with a simple coaction makes $U(\ggf)$ into a braided commutative algebra in the category of Yetter--Drinfeld modules over $\OO(G)$ and, therefore, the smash product algebra
	$\OO(G)\sharp U(\ggf)$ of regular differential operators carries the structure of a Hopf algebroid over (noncommutative) base algebras $\Ug^\mathrm{op}, U(\ggf)$.
	This construction retains essential features from a related construction of a Hopf algebroid structure on a noncommutative phase space of Lie type, considered in recent physics literature, while it avoids the need of completed tensor products and completions of algebraic structures. Steps of the construction are geometrically motivated, and are carried out for an affine algebraic group $G$ over any field. There a certain matrix $\OO$ of representative functions on $G$ is used. This matrix expresses an operator mapping the right invariant derivations to the corresponding left invariant derivations.  
	
	When $G$ is a Lie group over $\mathbb R$ or $\mathbb C$, analogous construction puts a Hopf algebroid structure over base algebras $\Ug^{\mathrm{op}}, \Ug$ on the subalgebra $\mathcal H \sharp \Ug$ of the algebra of differential operators, where $\mathcal H$ is some Hopf algebra of representative smooth functions on $G$. In this case the matrix $\OO$ has the same components as the matrix  of the adjoint representation $\Ad$ of the Lie group $G$ in a chosen basis. The same construction can be carried out also for any Hopf algebra $\mathcal H$ of germs of representative smooth functions around the unit of $G$. 
	
	We also give another related example of a Hopf algebroid structure over base algebras $\Ug^\op, \Ug$ on the finite dual Heisenberg double $\Ug^\circ \sharp \Ug$ for a Lie algebra $\ggf$ over any field. Here we use a certain matrix $\UU$ whose components are functionals derived from the adjoint representation $\ad$ of the Lie algebra $\ggf$. The constructions do not depend on the choice of  a basis of the Lie algebra $\ggf$.
	
	Moreover, the theorems hold for smash products $\mathcal H\sharp \Ug$ for any Hopf algebra $\mathcal H$ of representative functions on $G$ such that it contains the components of the matrix $\OO$, in the case with functions, and for any Hopf algebra $\mathcal H$ dual to $\Ug$ such that it contains the components of the matrix $\UU$, in the case with functionals. These algebras of the form $\mathcal H\sharp \UgL$ can also be presented as smash product algebras $\mathcal H^{\mathrm{co}} \sharp \UgR$, $\UgL \sharp \mathcal H^{\mathrm{co}}$  and $ \UgR \sharp \mathcal H$. We provide the formulas for the Hopf algebroid structure maps in each smash product algebra presentation.
\end{abstract}

\tableofcontents

\section{Introduction}
\subsection{Background on Hopf algebroids}
\subsubsection{Commutative Hopf algebroid}

If a group $(G,\cdot_G)$ acts from the right on a set $X$ one can form the corresponding action (or transformation) groupoid~\cite{ncmodels} with set of arrows $X\times G$ and set of objects $X$; if $G$ and $X$ are in some category of spaces, the groupoid will usually have the same kind of structure. Algebra $H = \mathrm{Fun}(G)$ of (appropriate class of) functions on $G$ has a structure of a Hopf algebra (possibly in a completed sense) $(H,\Delta,\epsilon)$ with comultiplication $f\mapsto\Delta(f) = \sum f_{(1)}\otimes f_{(2)}$, $\Delta(f)(a,b) = f(a\cdot_G b)$ (where some version of the formula $\mathrm{Fun}(G\times G)\cong\mathrm{Fun}(G)\otimes\mathrm{Fun}(G)$ is used) and $H$ acts on $A=\mathrm{Fun}(X)$ by a Hopf action (that is, $A$ is an $H$-module algebra, $h \btr(a\cdot b) = ( h_{(1)} \btr a)\cdot( h_{(2)} \btr b)$, $h \btr 1 =\epsilon(h) 1$). Applying the duality between function algebras and spaces, one obtains a commutative Hopf algebroid $\mathrm{Fun}(X\times G)\cong A\otimes H$ over $\mathrm{Fun}(X) = A$. 

\subsubsection{Hopf algebroid over noncommutative base} 

J-H.~Lu~\cite{Lu} introduced a left (associative) bialgebroid $(\mathcal{B},\alpha\colon A\to\mathcal{B},\beta\colon A^{\mathrm{op}}\to\mathcal{B},\Delta_{\mathcal{B}}\colon\mathcal{B}\to\mathcal{B}\otimes_A\mathcal{B},\epsilon_{\mathcal{B}}\colon\mathcal{B}\to A)$ over a noncommutative base algebra $A$, where $\mathcal{B}$ is an algebra, $\alpha,\beta$ algebra maps with mutually commuting images equipping $\mathcal{B}$ with a structure of an $A$-bimodule via $a.b.a' = \alpha(a)\beta(a')b$ and $\Delta_{\mathcal{B}}$ is a coassociative comultiplication in a category of $A$-bimodules with counit $\epsilon_{\mathcal{B}}$; these data are subject to standard axioms~\cite{Lu,bohmHbk,Bohm,BrzMilitaru,SS:two,xu}. Hopf algebroids are bialgebroids with an antipode antihomomorphism $\tau\colon\mathcal{B}\to\mathcal{B}$ which is tricky to axiomatize and few variants exist in the literature. Lu~\cite{Lu} also generalized commutative Hopf algebroids from action groupoids to Hopf algebroids over $A$ whose total algebra is smash product algebra $\mathcal{B} = A \sharp H$~\cite{bohmHbk,Lu,SS:two,stojic} with underlying vector space $A\otimes H$ and multiplication $(a\sharp h)\cdot(b\sharp k) = a (h_{(1)} \btr b) \sharp h_{(2)} k$, where $a\sharp h$ denotes $a\otimes h$ within the smash product. 

In a version of this construction from~\cite{BrzMilitaru}, $H$ is a Hopf algebra  and $A$ is a left-right braided commutative Yetter--Drinfeld $H$-module algebra. That means not only that the left action of $H$ on $A$ is Hopf, but there is a counital right coaction $\rho\colon A\to A\otimes H$ which is an algebra antihomomorphism such that two axioms, Yetter--Drinfeld property and braided commutativity, hold. This scalar extension bialgebroid is given by $\Delta_{\mathcal{B}}(a \sharp h) = a \sharp h_{(1)}\otimes 1\sharp h_{(2)}$, $\epsilon_{\mathcal{B}}(a \sharp h) = \epsilon(h)a$, $\alpha(a) = 1\sharp a$ and $\beta$ is $\rho$ followed by the identification of the underlying vector spaces of $A\otimes H$ and $A\sharp H$~\cite{bohmHbk,BrzMilitaru,stojic}.

\subsubsection{Heisenberg double}
If $A$ is a finite-dimensional Hopf algebra, the dual Hopf algebra $A^*$ and $A$ act on each other via Hopf actions which are canonically extended to Yetter--Drinfeld module algebra structures and the smash product $A\sharp A^*$ is called the Heisenberg double and serves as the main classical example of a scalar extension Hopf algebroid over $A$. In infinite-dimensional case, one usually has several candidates for a Hopf algebra $K$ in nondegenerate Hopf pairing with $A$; usually one needs to take $K$ as a topological Hopf algebra where the comultiplication takes values in some completion of $K\otimes K$. Neglecting issues of completions (and with expositional errors) article~\cite{heisd} has shown that the universal enveloping algebra $U(\ggf)$ together with corresponding deformed derivatives from physics literature generate a Heisenberg double of $U(\ggf)$; algebraic dual $H = U(\ggf)^*$ is a Hopf algebra in formally completed sense and $U(\ggf)$ has an $H$-coaction, so that it is (up to issues of completion) a Yetter--Drinfeld $H$-module algebra. Smash product $U(\ggf)\sharp H$ has been given a structure of a formally completed version of a $U(\ggf)$-Hopf algebroid (noncommutative phase space of Lie algebra type) in~\cite{halg}. A natural  version of a Heisenberg double of $U(\ggf)$ as an internal Hopf algebroid in a symmetric monoidal category $(\ipV, \tilde\otimes,k)$ of filtered cofiltered vector spaces has been fully developed in~\cite{stojicphd}. The canonical antiisomorphism between the Lie algebra $\ggf_L$ of left invariant and the Lie algebra $\ggf_R$ of right invariant vector fields can be viewed in coordinates as a matrix valued function $\OO$ and (due to the pairing between differential operators and smooth functions) identified with a matrix whose entries are in a dual of $U(\ggf_L)$ (or of $U(\ggf_R)$, if one wishes). 

\subsubsection{Adjoint representation} We observed  that the target
map $\beta$ restricted to $\ggf$ (in either version, \cite{halg} and \cite{stojicphd})
is obtained by multiplying the result of source map $\alpha$ with that matrix within the smash product. %(in general, the target map is often hardest to construct in new examples of bialgebroids)
This lead to the idea that one may take a smaller Hopf subalgebra within the dual of $U(\ggf_L)$ as long as the smash product contains also $U(\ggf_R)$ obtained by using $\beta$ and such that the comultiplication is taking values in non-completed tensor product. 
%Even more extreme, as the difference is only in the matrix of transition between the left and right invariant vector fields, one can try to model this {\em externally}, namely take an abstract transition matrix such that one gets the correct Lie algebra. The variety of possibilities is the automorphism group of the Lie algebra; hence there should be a bialgebroid that is a smash product of $U(\ggf)$ and $H =\mathcal{O}(\operatorname{Aut}(\ggf))$, where $\mathcal{O}(\operatorname{Aut}(\ggf))$ is a genuine (without completions) Hopf algebra of regular functions on the algebraic group $\operatorname{Aut}(\ggf)$ which makes sense over any field~$\genfd$. 
The formulas for the transition between the left and right invariant vector fields in a point is expressed with help of the adjoint action. It is not a surprise that the adjoint action could be used to construct an action type bialgebroid. Although the adjoint action is such a classical notion, we do not have the corresponding groupoid as only Hopf algebras give a common ground both for groups (via function Hopf algebras) and for Lie algebras (via universal enveloping algebras). Thus if an automorphism group $H$ (e.g.\ $\mathrm{Aut}(G)$) acts on a group $G = X$ (as a space) we obtain an action groupoid, but if we need to pass to infinitesimal side of the group space $X$, we may be able to see the action information at the level of bialgebroids.

\subsection{Our results}
\subsubsection{For an affine algebraic group over any field}
Let $G$ be an affine algebraic group over any field~$\genfd$ and denote by $\OO(G)$ the algebra of regular functions on $G$. Denote by $\gL$ and $\gR$ the Lie algebra of left invariant derivations of $\OO(G)$ and the Lie algebra of right invariant derivations of $\OO(G)$ respectively. We prove that $U(\gL)$ is a right-left braided commutative Yetter--Drinfeld module algebra over $\OO(G)$. Here the right action $\btl$ of $f\in\OO(G)$ on $D\in\UgL$ is defined to satisfy $(D \btl f)(k) = D(fk)$ for all $k\in \OO(G)$ and the left coaction is geometrically defined such that it maps left invariant derivations into the corresponding right invariant ones written as elements of the smash product algebra $\OO(G)\sharp \UgL$ induced by the above action. This uses a certain matrix $\OO$ whose components are the components of the adjoint map $\Ad$ and puts a structure of a Hopf algebroid over base algebras $\UgR, \UgL$ on the smash product algebra $\OO(G)\sharp \UgL$. Likewise, the right action of $\OO(G)^\co$ on $\UgR$ defined by the same  formula and the left coaction that maps right invariant derivations to the corresponding left invariant ones as elements of $\OO(G)^\co \sharp \UgR$ make $\UgR$ into a braided commutative Yetter--Drinfeld module algebra over $\OO(G)^\co$. Both algebras $\OO(G)\sharp \UgL $ and $\OO(G)^\co \sharp \UgR$ are isomorphic to the algebra $\operatorname{Diff}(G)$ of regular differential operators on the affine algebraic group $G$. Therefore, by results in \cite{stojic}, $\operatorname{Diff}(G)$ is a scalar extension Hopf algebroid over base algebras $\UgR,\UgL$ and in a different way, also over base algebras $\UgL,\UgR$. Moreover, we prove that the same is true for $\mathcal H\sharp \UgL$ and $\mathcal H^\co \sharp \UgR$ for any Hopf algebra $\mathcal H$ such that $\OO^{min}(G) \subset \mathcal H \subset \OO(G)$, where $\OO^{min}(G) \subset \OO(G)$ is the smallest Hopf subalgebra containing the components of the matrix of the adjoint map $\Ad$ with regard to a basis of $\ggf$. The constructions do not depend on the choice of the basis of the Lie algebra $\ggf$.
\subsubsection{For a Lie group}
The analogous claims are proven for a Lie algebra $G$ over $\mathbb R$ or $\mathbb C$, and $\mathcal H$ some Hopf algebra of representative smooth functions on $G$.
\subsubsection{For the finite dual of a universal enveloping algebra} 
We also present related examples of scalar extensions of $\Ug$ with subalgebras of $\Ug^*$ for any finite-dimensional Lie algebra $\ggf$ over any field $\genfd$. We prove that $U(\ggf)$ is a braided commutative right-left Yetter--Drinfeld module algebra over the finite dual $\Ug^\circ$. Here the right action $\btl$ of $f\in\Ug^\circ$ on $D\in\Ug$ is defined to satisfy $\langle D \btl f, k \rangle = \langle D, fk \rangle$ for all $k\in \Ug^\circ$ and the left coaction is defined by introducing a certain matrix $\UU$ of functionals, with the components defined by formulas involving the matrices of the adjoint representation of $\ggf$ on $\ggf$, with properties analogous to the properties of the matrix $\OO$ above. This puts a structure of a Hopf algebroid over base algebras $\Ug^\op, \Ug$ on the smash product algebra $\Ug^\circ\sharp \Ug$. Moreover, we prove that the same is true for $\mathcal H\sharp \Ug$ for any Hopf algebra $\mathcal H$ such that $\Ug^{min} \subset \mathcal H \subset \Ug^\circ$, where $\Ug^{min}$ is the Hopf algebra generated by the components of the matrix $\UU$. The constructions do not depend on the choice of the basis of the Lie algebra $\ggf$.

\section{Preliminaries}
Throughout the paper we freely use Sweedler notation and the Kronecker symbol $\delta^i_j$.  Symbol $\otimes$ denotes the tensor product of vector spaces over the ground field $\genfd$.

\subsection{Yetter--Drinfeld module algebra}

\begin{definition}\label{rlydma} 
	Let $H$ be a $\genfd$-bialgebra. A \emph{right-left Yetter--Drinfeld module algebra} over $H$ is a triple $(A,\btl,\lambda)$ such that 
	\begin{enumerate} 
		\item[(1)] $A$ is a $\genfd$-algebra 
		\item[(2)] $(A,\btl)$ is a right $H$-module  
		\item[(3)] $(A,\lambda\colon a\mapsto a_{[-1]}\otimes a_{[0]})$ is a left $H$-comodule		
		\item[(4)] action and coaction satisfy the \emph{right-left Yetter--Drinfeld condition} 
		\begin{equation} \label{eq:YDrl1}
		\sum a_{[-1]} h_{(1)}\otimes (a_{[0]} \btl h_{(2)}) = \sum h_{(2)}(a \btl h_{(1)} )_{[-1]} \otimes  (a \btl h_{(1)})_{[0]}, \quad \forall h \in H, a\in A
		\end{equation}	
		\item[(5)] action $\btl$ is \emph{Hopf}, that is 
		\begin{equation}
		(ab) \btl h = \sum 	(a \btl h_{(1)})(b \btl h_{(2)}), \quad 1_A \btl h = \epsilon(h) 1_A, \quad \forall h\in H, a,b\in A
		\end{equation}		
		\item[(6)] coaction $\lambda$ satisfies 
		\begin{equation} \label{eq:marl}
		\sum (ab)_{[-1]} \otimes (ab)_{[0]} = \sum b_{[-1]}a_{[-1]} \otimes a_{[0]} b_{[0]}, \quad \lambda(1_A) = 1_H\otimes 1_A, \quad \forall a,b\in A
		\end{equation} 
		called the \emph{comodule algebra property} (over $H^\op$).
	\end{enumerate}
	We say that a Yetter--Drinfeld module algebra $(A,\btl,\lambda)$ is \emph{braided commutative} if
	\begin{equation} \label{eq:bcrl}
	\sum (a \btl b_{[-1]})b_{[0]} = ba, \quad \forall a,b\in A.
	\end{equation} 
\end{definition}

Axioms (2), (3) and (4) from Definition \ref{rlydma} together comprise the axioms for $(A,\btl,\lambda)$ being a \emph{right-left Yetter--Drinfeld module} over $H$. Morphisms of Yetter--Drinfeld $H$-modules are $H$-module morphisms which are also $H$-comodule morphisms. Right-left Yetter--Drinfeld $H$-modules form a category ${}^H\mathcal{Y}\mathcal{D}_H$ carrying a canonical monoidal structure with a pre-braiding, which is a braiding if $H$ is a Hopf algebra. 

Axioms (1), (2) and (5) from Definition \ref{rlydma} together say that $(A,\btl)$ is a \emph{right module algebra} over $H$.  Given any right module $H$-algebra $(A,\btl)$, the vector space $H\otimes A$ carries a structure of an associative $\genfd$-algebra with the multiplication bilinearly extending formula
\begin{equation}
(h\otimes a)\cdot(h'\otimes a') = \sum h h'_{(1)}\otimes (a\blacktriangleleft h'_{(2)}) a', \quad  \text{ for } h,h'\in H \text{ and } a,a' \in A
\end{equation} 
and the unit $1_H\otimes 1_A$. This algebra is called the \emph{smash product algebra} and is denoted	$H\sharp_{\btl} A$. If it is clear from the context which action induces the smash product multiplication, we denote this smash product algebra simply by $H\sharp A$. Smash product algebra $H\sharp A$ comes along with the canonical algebra monomorphisms $A\hookrightarrow H \sharp A$ and $H\hookrightarrow H\sharp A$. The images of these two embeddings are denoted $1\sharp A$ and $H\sharp 1$. A general element $a\otimes h$ of $A\sharp H$ is denoted by $a\sharp h$.

Written in terms of the corresponding smash product algebra $H\sharp A$, the right-left Yetter--Drinfeld property is
\begin{align}
\lambda(a) \cdot h
&= h_{(2)} \cdot \lambda(a \btl h_{(1)} ), \quad  \text{ for all } a\in A \text{ and } h\in H. \label{eq:YDrl}
\end{align}

Axioms (1), (3) and (4) from Definition \ref{rlydma} define a \emph{left $H^\op$-comodule algebra} $(A,\lambda)$. Therefore, a right-left Yetter--Drinfeld $H$-module algebra is a right-left Yetter--Drin\-feld $H$-module and an algebra that is a right $H$-module algebra and a left $H^\op$-comodule algebra.

\subsection{Hopf pairing}
Pairings of vector spaces are bilinear maps into the ground field $\genfd$ which are in this work not required to be nondegenerate. 
\begin{definition} 
	A pairing between two Hopf algebras $A$ and $H$ is a \emph{Hopf pairing} if
	$\langle \Delta_A(a), h\otimes k\rangle = \langle a, h\cdot k\rangle$ where on the left hand side the pairing is a product of pairings in each tensor factor and likewise $\langle a\otimes b, \Delta_H(h)\rangle = \langle a\cdot b,h\rangle$,
	as well as $\epsilon_A(a) = \langle a,1\rangle$, $\epsilon_H(h) = \langle 1,h\rangle$ and $\langle S_A(a), h \rangle = \langle a, S_H(h)\rangle$, for all $a\in A$ and $h,k\in H$.
\end{definition}

If $A$ is a Hopf $\genfd$-algebra and $(H,\str)$ is a left module $A$-algebra, the vector space $H\otimes A$ carries a structure of an associative $\genfd$-algebra with the multiplication bilinearly extending formula
\begin{equation}
(h\otimes a)\cdot( h'\otimes a') = \sum h (a_{(1)}\str h')\otimes a_{(2)} a', \quad \text{ for } h,h'\in H \text{ and } a,a' \in A
\end{equation} 
and the unit $1_H\otimes 1_A$. This algebra is also called the \emph{smash product algebra} and is denoted by $H\sharp_{\str} A$, or simply $H\sharp A$ if the action is clear from the context. 

\begin{proposition}\label{prop:leftright}
	A Hopf pairing $\langle \,,\rangle \colon A \otimes H \to \genfd$ between Hopf algebras $A$ and $H$ induces a left action and a right action
	\begin{align}
	\str &\colon A\otimes H \to H, \quad a\str h = \sum h_{(1)}\langle a, h_{(2)} \rangle , \\
	\btl &\colon A \otimes H \to A, \quad a \btl h = \sum \langle a_{(1)}, h \rangle a_{(2)}, \quad \text{for } a\in A \text{ and } h\in H,
		\end{align} which both induce the same smash product algebra $H\sharp A$. 
\end{proposition}
\begin{proof}
	The proof is straightforward. For $a,b \in A$ and $f,h\in H$, we have that 
	\begin{align*}
	f \sharp a \cdot h \sharp b & = \sum f a_{(1)} \langle a_{(1)}, h_{(2)} \rangle a_{(2)}b \\
	& = \sum f (a_{(1)} \str h) \sharp a_{(2)} b \\
	& = \sum f h_{(1)} \sharp (a \btl h_{(2)})b.
	\end{align*}

\end{proof}
\subsection{Hopf algebroid}

Given a Hopf algebra $H$, each right-left Yetter-Drinfeld $H$-module algebra $A$ gives rise to an associative $A$-bialgebroid structure on the smash product algebra $H \sharp A$ \cite{BrzMilitaru, stojic}. There is also an antimultiplicative antipode map making $H\sharp A$ the total algebra of a Hopf algebroid over $A$ \cite{BrzMilitaru,stojic}. A generalized version of the scalar extension construction is provided in \cite{stojic}. Its input is a compatible pair $A$ and $A^{\op}$ of a right-left and a left-right braided commutative Yetter-Drinfeld $H$-module algebras, and the output is a symmetric Hopf algebroid $A^\op\sharp H \cong H\sharp A$ over base algebras $A^\op,A$. The construction in \cite{stojic} does not require that the antipode of $H$ is invertible. 

\subsection{Representative functions} 
If $G$ is a group, the multiplication $m$ of $G$ induces a map $m^* \colon \mathcal F(G,\genfd) \to \mathcal F(G\times G,\genfd)$ by $m^*(f) = f \circ m$, for a function $f \in \mathcal F(G,\genfd)$. We say that $f \in \mathcal F(G, \genfd)$ is a \emph{representative function} if $m^*(f)$ is in the image of the canonical map $\mathcal F(G,\genfd)\otimes \mathcal F(G,\genfd) \to \mathcal F(G\times G,\genfd)$ \cite{HochMost}. The set $\mathcal R$  of representative functions on a group $G$ is a Hopf algebra with the structure maps for the comultiplication, the counit and the antipode transposed from the structure maps of $G$, 
$$\Delta(f)(a\otimes b) = f(ab), \quad \epsilon(f)=f(e), \quad S(f)(a)= f(a^{-1}) ,\quad \text{ for } f\in \mathcal R \text{ and } a,b\in G .$$
If $G$ is an affine algebraic group over a field $\genfd$, then $\genfd$-valued representative functions on $G$ are automatically regular. If $G$ is a Hausdorff topological group, then $\genfd$-valued representative functions on $G$ are automatically continuous.   %Equivalenty, $f$ is representative if the span of all functions $\{g \cdot f : g\in G\}$, where $g \cdot f$ denotes the function with mapping rule $a \mapsto f(ag)$ for $a \in A$, is finite dimensional. 

\section{Overview of results}\label{sec:cases}

We consider the following cases:
\begin{enumerate}
	\item[(I)] $G$ is a Lie group over a field $\genfd$ equal to $\mathbb{R}$ or $\mathbb{C}$, $\mathcal{L}$ is the tangent space at the unit of $G$, and
	\begin{enumerate} 
		\item[(a)] $\mathcal{F}$ is the algebra $C^\infty(G)$ of smooth functions on $G$, $\mathcal{D}$ is the Lie algebra %$\Gamma^\infty(G)$ 
		of smooth vector fields on $G$,   or
		\item[(b)] $\mathcal{F}$ is the algebra $C^\infty(G,e)$ of germs of smooth functions around the unit of $G$, $\mathcal{D}$ is the Lie algebra of germs of smooth vector fields around the unit on $G$, or
	\end{enumerate}
	\item[(II)] $G$ is an affine algebraic group over any field $\genfd$, $\mathcal{F}$ is the Hopf algebra $\OO(G)$ of regular functions on $G$, $\mathcal{L}$ is the algebra of differentiations $\OO(G) \to \genfd$, $\mathcal{D}$ is the Lie algebra of derivations $\OO(G) \to \OO(G)$, or
	\item[(III)] \begin{enumerate}
		\item[(a)] $G$ is a Lie group over a field $\genfd$ equal to $\mathbb{R}$ or $\mathbb{C}$, $\mathcal{L}$ is the tangent space at the unit of $G$, $\mathcal{F}$ is the algebra $J^\infty(G,e)$ of $\infty$-jets of smooth functions around the unit of $G$,  $\mathcal{D}$ is the Lie algebra  of $\infty$-jets  of smooth vector fields around the unit of $G$, or
		\item[(b)] $\mathfrak{g}$ is a Lie algebra over any field $\genfd$, $\mathcal{F}$ is the dual $\Ug^*$ of the universal enveloping algebra $U(\ggf)$.
	\end{enumerate}
\end{enumerate} 
In cases (I) and (II), we consider algebras $\mathcal H \subset \mathcal F \cap \mathcal R$, that is, Hopf algebras of some representative functions on $G$, and in case (III) $\mathcal H$ is a Hopf subalgebra of $\Ug^\circ$.
We prove that $\Ug$ is a braided commutative right-left Yetter--Drinfeld module algebra over $\mathcal H$ if $\mathcal H$ contains as elements components of a certain matrix $\OO$ in cases (I) and~(II), or components of a certain matrix $\UU$ in case (III). Matrix $\OO$ is a
the matrix of the adjoint representation $\Ad$ of group~$G$, and matrix $\UU$ is likewise derived from adjoint representation $\ad$ of Lie algebra $\ggf$, with regard to a chosen basis of $\ggf$, and the resulting minimal Hopf algebras $\OO^{min}(G)$ and $\Ug^{min}$ that contain their components do not depend on the choice of the basis for $\ggf$. This makes $\mathcal H \sharp \Ug \cong \Ug^\op \sharp \mathcal H$ into a scalar extension Hopf algebroid over base algebras $\Ug^\op, \Ug$ \cite{stojic}.

\section{Yetter--Drinfeld module algebra over a Hopf algebra of functions}

We consider now cases (I) and (II) described in Section \ref{sec:cases}. In case (II), $\mathcal F$ is a Hopf algebra. In case (I) we usually work with some Hopf algebra $\mathcal H \subset \mathcal F$ of representative functions, and in case (II) with whole $\mathcal F$.

\subsection{Hopf algebra of some representative functions on $G$}
Let $\mathcal{H}$ be a subalgebra of $\mathcal{F}$ and a Hopf algebra such that its coalgebra structure and antipode are the \emph{transpose of the group structure} of $G$, that is $\Delta_\mathcal{H}(f)(a\otimes b) = f(ab)$, $\epsilon_\mathcal{H}(f) = f(e)$, and $S_\mathcal{H}(f)(a) = f(a^{-1})$, for all $f\in \mathcal{H}$ and $a,b \in G$. In other words, $\mathcal H$ is a subset of the Hopf algebra $\mathcal R \cap \mathcal F$ of representative functions on $G$, with the inherited Hopf algebra structure. Note that even though the algebra $\mathcal F$ in case (I) of smooth functions is not a Hopf algebra, maps $\epsilon$ and $S$ are defined, as follows. Denote by $\epsilon$ the $\genfd$-linear map $\mathcal F \to \genfd$ such that $\epsilon(f) = f(e)$ for all $f\in\mathcal F$, and by $S$  the $\genfd$-linear map $\mathcal F \to \mathcal F $ such that $S(f)(a) = f (a^{-1}) $ for all $f\in\mathcal F$ and $a\in G$. In case (II), $\mathcal H$ is any Hopf subalgebra of $\mathcal F$ since $\mathcal F$ is a subset of the Hopf algebra $\mathcal R$ of representative functions on $G$.

\subsection{Left invariant and right invariant derivations} 
A \emph{derivation} of a $\genfd$-algebra $\mathcal F$ is any $\genfd$-linear map $D \colon \mathcal F \to \mathcal F$ such that $$D(fh) = D(f)h + f D(h), \quad \text{ for all } f,h \in \mathcal F.$$ Denote by $\mathcal D$ the Lie algebra of derivations of $\mathcal F$. The Lie algebra bracket is simply the commutator of elements. 
Let $\gL$ denote the Lie algebra of \emph{left invariant} derivations, that is, elements $D \in \mathcal{D}$ such that for every $g\in G$ we have that $D(f\circ L_g) = D(f) \circ L_g$, where $L_g \colon G \to G $, $L_g(a) = ga$ for $a \in G$. Let $\gR$ denote the Lie algebra of \emph{right invariant} elements of $\mathcal{D}$, that is, elements $D \in \mathcal{D}$ such that for every $g\in G$ we have that $D(f\circ R_g) = D(f) \circ R_g$, where $R_g \colon G \to G $, $R_g(a) = ag$ for $a \in G$. It is easy to see that $\gL$ and $\gR$ are indeed  Lie algebras as subalgebras of the Lie algebra $\mathcal D$. Denote by $\UgL$ and $\UgR$ the universal enveloping algebras of $\gL$ and $\gR$ respectively. 

Denote by $\mathcal L$ the $\genfd$-vector space of differentiations of $\mathcal F$. A \emph{differentiation} of $\mathcal F$ is any $\genfd$-linear map $\nu \colon \mathcal F \to \genfd$ such that $\nu(fg) = \nu(f) \epsilon(g) + \epsilon(f) \nu(g)$. All differentiations on $\mathcal F$ comprise a "tangent space" $\mathcal{L}$ at the unit of $G$. It is finite-dimensional in the case of an affine algebraic group $G$ \cite{Hochschild}. Result $\nu(f)$ is often denoted $\langle \nu, f \rangle$, for $\nu \in \mathcal L$ and $f \in \mathcal F$, since $\nu$ is a linear functional on $\mathcal F$.

\begin{proposition}
	
	Dimension of $\genfd$-vector space $\mathcal{L} $ of differentiations of $\mathcal F$ is finite.
\end{proposition}
\begin{proof}
	For case (II), vector $\genfd$-space $\mathcal L$ is finite-dimensional, because a differentiation of $\mathcal F$ is determined by its values on a finite system of $\genfd$-algebra generators \cite{Hochschild}.  For case (I), this is the tangent space at the unit of a Lie group. 
\end{proof}

If $D\in \mathcal D$, then $\epsilon \circ D \in \mathcal L$. For case~(I), differentiation $\epsilon\circ D$ is the tangent vector at the unit of $G$ that corresponds to the vector field $D$. Therefore we have a $\genfd$-linear map $\epsilon \circ - \colon \mathcal D \to \mathcal L$. The next proposition proves that this map is surjective.

\begin{proposition}
For each differentiation $\nu \in \mathcal L$ there exists a unique left invariant derivation $X_\nu$ and a unique right invariant derivation $Y_\nu$ such that $\epsilon \circ X_\nu = \epsilon \circ Y_\nu = \nu$. In case (II), they are given by formulas: $$X_\nu(f) = \sum  f_{(1)} \langle \nu, f_{(2)}  \rangle, \quad Y_\nu(f) = \sum f_{(1)} \langle \nu, f_{(2)} \rangle, \text{ for all } f\in \mathcal F.$$ In case (I), the above formulas are true for elements $f$ of the Hopf algebra $\mathcal H \subset \mathcal F \cap \mathcal R$.
\end{proposition}
\begin{proof} The right invariant derivation $Y_\nu$ is uniquely determined by the formula $Y_\nu(f)(g) = \nu(f\circ R_g)$ for all $f\in \mathcal F$ and $g\in G$, and similarly for the left invariant one. For case (I), $Y_\nu$ is the right invariant vector field with the corresponding tangent vector $\nu$ at the unit, and similarly for $X_\nu$. For any Hopf algebra $\mathcal H \subset \mathcal F$, in both cases (I) and (II),  for $f \in \mathcal H$, $Y_\nu(f)(a) = Y_{\nu}(f)(e \cdot a) = Y_\nu(f \circ R_a)(e) = \sum \langle \nu, f_{(1)} f_{(2)}(a) \rangle = \sum \langle \nu, f_{(1)} \rangle f_{(2)}(a),$ and therefore $Y_\nu(f) = \sum \langle \nu, f_{(1)} \rangle f_{(2)}$, and similarly for $X_\nu$  \cite{Hochschild}.
\end{proof}

Lie algebra $\mathcal D$ of derivations is a left $\mathcal F$-module in a usual way.  
\begin{proposition} \label{prop:isoleft} The mapping  $$\mathcal{F} \otimes \mathcal{L} \to \mathcal{D}, \quad f \otimes \nu \mapsto fX_\nu,$$ is an $\mathcal F$-module isomorphism.  The image of $\genfd \otimes \mathcal L$ is $\gL \subset \mathcal{D}$.
\end{proposition}
\begin{proof}
	See \cite{Hochschild}, Theorem III.3.1.
\end{proof}	
\begin{proposition} \label{prop:isoright} The mapping  $\mathcal{F} \otimes \mathcal{L} \to \mathcal{D}$ defined by the rule $f \otimes \nu \mapsto fY_\nu$ is an $\mathcal F$-module isomorphism.  The image of $\genfd \otimes \mathcal L$ is $\gR \subset \mathcal{D}$.
\end{proposition}
\begin{proof}
	Similar as Proposition \ref{prop:isoleft}.
\end{proof}
It follows that each right invariant derivation can be written as a product of functions and elements of a basis of the Lie algebra of left invariant derivations in a unique way, and likewise for left invariant derivations. It can be proven that functions that appear thereat as coefficients are always representative.

\subsection{Hopf actions and pairings}

\begin{proposition}\label{prop:action}
	The action $\mathcal D \otimes \mathcal F \to \mathcal F$ restricted to $\gL \otimes \mathcal F \to \mathcal F$ induces a unique left Hopf action $$\str \colon \UgL \otimes \mathcal F \to \mathcal F,$$ that is, differentiation of a function by a left invariant differential operator. 

\end{proposition}
\begin{proof} 
	In case (I), every left invariant vector field  $X \in \ggf_L$ defines a first order differential operator $D_X \colon C^\infty(G) \to C^\infty(G)$ by $D_Xf = Xf$. By multiplicativity and linearity, the mapping $X \mapsto D_X$ extends to an injective homomorphism of unital algebras $$D \colon \UgL \hookrightarrow \End_{\mathbb{R}}(C^\infty(G)), \quad P \mapsto D_P.$$ Explicitly, element  $P = X_1\cdots X_r \in U(\ggf_L)$, where $X_1,\ldots,X_r \in \ggf_L$, maps to a diferential operator $D_{P}$ defined by  %\ref{Timmermann}
$$(D_{P}f)(y) = (D_{X_1} \cdots D_{X_r}f)(y) 
= \lim_{t_1,\ldots,t_r\to 0}\frac{\partial^{r}}{\partial t_1 \cdots \partial t_r} f(y \cdot \exp(t_1X_1)\cdots\exp(t_rX_r)).
$$
This formula is obtained by starting with the formula for a left invariant vector field~$X$,
$$(D_X f)(y) = \lim_{t\to 0}\frac{\partial}{\partial t} f(y \cdot \exp(t X)),$$ 
where the argument $y$ is multiplied on the right with the exponential for the formula to be left invariant, and by using induction on the order of the differential operator, with the step of induction simply
$$(D_{X_1} \cdots D_{X_r}f)(y) 
= \lim_{t\to 0}\frac{\partial }{\partial t_1} (D_{X_2} \cdots D_{X_r}f)(y \cdot \exp(t_1X_1)) .$$
%The aforementioned pairing is defined as	$$\langle\,,\rangle \colon U(g_L) \otimes \mathcal F \to \genfd, \quad \langle P,f\rangle = (D_P f)(e),$$
%the pairing of Hopf algebra  $U(\ggf_L)$ and algebra $\mathcal F$. 
This action is Hopf due to Leibniz rule.  

In case (II), we define the action of $T(\gL)$ in the same way, simply by composing left invariant derivations. Let $X_1,\ldots,X_n$ be a basis of $\gL$. The derivation $X_iX_j - X_j X_i$ is again a left invariant derivation and acts in the same way as $[X_i,X_j]$ because it is the same. Hence we have a well defined action of $\UgL$ on $\mathcal F$. It is easy to check that the action is Hopf.
\end{proof}

\begin{proposition}\label{prop:actionR}
	The action $\mathcal D \otimes \mathcal F \to \mathcal F$ restricted to $\gR \otimes \mathcal F \to \mathcal F$ induces a unique right Hopf action $$\str \colon \UgR \otimes \mathcal F \to \mathcal F,$$ that is, differentiation of a function by a right invariant differential operator. 
	
\end{proposition}
\begin{proof}
	Analogous to the proof of Proposition \ref{prop:action}. This pairing is defined by using formula
	$$(D_{P}f)(y) = (D_{Y_1} \cdots D_{Y_r}f)(y) 
	= \lim_{{t_1,\ldots,t_r\to 0}}\frac{\partial^{r}}{\partial t_1 \cdots \partial t_r} f(\exp(t_rY_r)\cdots\exp(t_1Y_1) \cdot y)
	$$
	for $P = Y_1\cdots Y_r \in U(\ggf_R)$, where $Y_1,\ldots,Y_r \in \ggf_R$. Note that the exponentials here are on the left and that the order of the arguments is reversed, due to $P$ being right invariant instead of left invariant.
\end{proof}

The pairing of the universal enveloping algebra and the algebra of functions is defined as applying the invariant differential operator to a function and then evaluating the result at the unit.

\begin{proposition}\label{prop:pair} The Hopf action $\str \colon \UgL \otimes \mathcal F \to \mathcal F$ from Proposition \ref{prop:action}
	induces a pairing 
	$$\langle\,,\rangle \colon \UgL \otimes \mathcal{F} \to \genfd,$$ 
	by formula
	\begin{equation}\label{eq:pair}  \langle D, f\rangle = \epsilon(D(f)), \quad \text{ for all } D \in \UgL \text{ and } f \in \mathcal F.
	\end{equation}  
	Its restriction  to a pairing with a Hopf algebra $\mathcal H \subset \mathcal F \cap \mathcal R$ is a Hopf pairing $$\langle\,,\rangle \colon \UgL \otimes \mathcal{H} \to \genfd.$$
	Moreover, action $\str$ satisfies the formula
	\begin{equation}\label{eq:actpair} D \str h = \sum h_{(1)} \langle D,h_{(2)} \rangle, \quad \text{ for all } D\in \UgL \text{ and } h \in \mathcal{H}
	\end{equation}
	and it restricts to a left Hopf action  
	$\str \colon \UgL \otimes \mathcal H \to \mathcal H.$
\end{proposition}

\begin{proof}

	For case (I), first we check the compatibility of the multiplication on $\UgL$ and the comultiplication on $H$,
\begin{align*}
\langle {X_1}\cdots{X_r}, f \rangle & = \lim_{t_1,\ldots,t_r\to 0}\frac{\partial^{r}}{\partial t_1 \cdots \partial t_r} f(\exp(t_1X_1)\cdots\exp(t_rX_r)) 
\\
& = \lim_{t_1,\ldots,t_r\to 0} \frac{\partial^{r}}{\partial t_1 \cdots \partial t_r}\left(\sum f_{(1)}(\exp(t_1X_1))\cdots f_{(r)}(\exp(t_rX_r))\right) 
\\
&  = \sum \frac{\partial}{\partial t_1 }\Big|_{t_1=0}  f_{(1)}(\exp(t_1X_1)) \cdots \frac{\partial}{\partial t_r}\Big|_{t_r=0} f_{(s)}(\exp(t_rX_r)) 
\\
& = \sum \langle X_1, f_{(1)}\rangle\cdots \langle X_r , f_{(r)} \rangle.
\end{align*}
Second, we check the compatibility of the multiplication on $H$ and the comultiplication on  $\UgL$. For $X\in \ggf_L$ we have Leibniz rule
$$ \langle \Delta(X), fg \rangle = f(e) \langle X,g\rangle + \langle X, f \rangle g(e)
$$
and it is easy to show inductively that the claim  holds for all  $P \in \UgL$.
Third, we check the compatibility of units and counits, for $P = X_1 \cdots X_r$:
$$\langle P, 1_H \rangle = 0 = \epsilon(P), \quad \langle 1_{\Ug}, f\rangle = f(e) = \epsilon(f)  $$
$$\langle 1_{\Ug}, 1_\mathcal H \rangle = 1 = \epsilon(1_{\Ug}).$$
Finally, we check the compatibility of antipodes: for $X \in \ggf_L$  we have
$$\langle X, Sf \rangle = \lim_{t\to 0} \frac{\partial}{\partial t }  f((\exp(tX))^{-1}) =\lim_{t\to 0} \frac{\partial}{\partial t }   f(\exp(-tX)) = - \langle X,f\rangle = \langle SX, f\rangle$$
and it is easy to prove inductively that the claim holds for all  $P \in \UgL$.

For case (I), we prove the formula, for $P \in \UgL$ and $f \in H$,
\begin{equation}
D_P f = \sum f_{(1)} \langle P, f_{(2)}\rangle.
\end{equation}

We check this first for $X \in \ggf_L$,
\begin{align*}
(D_X f)(y) &= \lim_{t\to 0} \frac{\partial}{\partial t} f(y \cdot \exp(tX)) \\
& = \lim_{t\to 0} \frac{\partial}{\partial t} \sum f_{(1)}(y) f_{(2)} (\exp(tX)) \\
& = \sum f_{(1)}(y) \langle X, f_{(2)} \rangle
\end{align*}
and then the claim is easily proven inductively for all $P \in \UgL$. 

For case (II), define the pairing by formula (\ref{eq:pair}). It is straightforward to check similarly as above that it is Hopf pairing and that it satisfies formula (\ref{eq:actpair}). 
\end{proof}

Since $\mathcal H$ is commutative, $\mathcal{H}^\co$ is again a Hopf algebra over $\genfd$. 

\begin{proposition} Hopf action $\str \colon \UgR \otimes \mathcal F \to \mathcal F$ from Proposition \ref{prop:actionR} induces a pairing 
	$$\langle\,,\rangle \colon \UgR \otimes \mathcal{F} \to \genfd$$ by formula 
	\begin{equation} \langle D, f\rangle = \epsilon(D(f)), \quad \text{ for all } D \in \UgR \text{ and } f \in \mathcal F.
	\end{equation} %The pairing thus consists of applying a right invariant differential operator to a function and then evaluating the result at the unit.
	Its restriction  to a pairing with a Hopf algebra $\mathcal H \subset \mathcal F \cap \mathcal R$ is a Hopf pairing $$\langle\,,\rangle \colon \UgR \otimes \mathcal{H}^\co \to \genfd.$$ 
	Moreover, action $\str$ satisfies the formula
	\begin{equation}
	D \str h = \sum \langle D,h_{(1)} \rangle h_{(2)} , \quad \text{ for all } D\in \UgR \text{ and } h \in \mathcal{H}
	\end{equation}
	and it restricts to a right Hopf action  
	$\str \colon \UgR \otimes \mathcal H \to \mathcal H.$

\end{proposition}
\begin{proof}
	
Similar as proof of Proposition \ref{prop:pair}.
\end{proof}

\subsection{Algebra of differential operators}

\begin{proposition}\label{prop:smash} The actions  $\str \colon \UgL \otimes \mathcal F \to \mathcal F$ and $\str \colon \UgR \otimes \mathcal F \to \mathcal F$ induce smash product algebra structures $\mathcal F \sharp_{\str} \UgL$ and $\mathcal F \sharp_{\str} \UgR$ on the corresponding tensor products. As algebras
	$$\operatorname{Diff}(G) \cong \mathcal F \sharp_{\str} \UgL \cong \mathcal F \sharp_{\str} \UgR.$$
Both actions extend to actions of the whole smash product algebra such that $\mathcal F\sharp 1$ acts on $\mathcal F$ by multiplication. These actions are precisely actions of differential operators on functions.

\end{proposition}
\begin{proof}
The multiplication in the smash product algebra is by definition $$f \sharp P \cdot h\sharp D = \sum f (P_{(1)} \str h) \sharp P_{(2)} D, \quad \text{ for } f,h\in\mathcal F \text{ and } P,D \in \UgL \text{ or } \UgR.$$ 	The proof is straightforward. 
\end{proof}

\subsection{Definition of Hopf algebra $\OO^{min}(G)$}

\begin{proposition}\label{minus}
	Let $X_1,\ldots, X_n$ be a basis of the Lie algebra $\gL$  of left invariant derivations of $\mathcal F$, and denote by $Y_1,\ldots, Y_n$ the corresponding basis of the Lie algebra $\gR$  of right invariant derivations of $\mathcal F$, that is, the one such that $\epsilon \circ X_i = \epsilon \circ Y_i$ for every $i \in \{1,\ldots,n\}$. Define the structure constants $C_{ij}^k$ from $$[X_i,X_j] = \sum_k C_{ij}^k X_k, \quad \text{ for } i,j\in\{1,\ldots,n\}.$$ Then $$[Y_i,Y_j] = - \sum_k C^k_{ij} Y_k, \quad \text{ for all } i,j\in\{ 1,\ldots,n\}.$$
	
\end{proposition}
\begin{proof} For case (I), this is well known. For case (II), $\mathcal F$ is a Hopf algebra and for every $f \in \mathcal F$ and every $i,j\in\{1,\ldots,n\}$, we have  that 
	\begin{align*}
	[Y_i,Y_j](f) & = \sum \langle Y_iY_j - Y_jY_i, f_{(1)} \rangle f_{(2)} \\
	& = \sum \langle Y_i, f_{(2)}\rangle \langle Y_j, f_{(1)} \rangle f_{(3)} - \sum \langle Y_j, f_{(2)} \rangle \langle Y_i, f_{(1)}\rangle f_{(3)} \\
	&= \sum \langle X_i, f_{(2)}\rangle \langle X_j, f_{(1)} \rangle f_{(3)} - \sum \langle X_j, f_{(2)} \rangle \langle X_i, f_{(1)}\rangle f_{(3)} \\
	&= \sum \langle X_jX_i - X_i X_j,f_{(1)}\rangle f_{(2)} = \sum \langle C_{ji}^k X_k,f_{(1)}\rangle f_{(2)} \\
	& = \sum C_{ji}^k \langle Y_k, f_{(1)}\rangle f_{(2)} = \left(\sum C_{ji}^k Y_k\right)(f).
	\end{align*}
\end{proof}
\begin{proposition}\label{prop:comm}
	Every left invariant derivation commutes with every right invariant derivation  in $\mathcal{D}$.	
	
\end{proposition} 
\begin{proof}
	For case (I), this is well known. For case (II), see Proposition I.2.2. in \cite{Hochschild}. 	
\end{proof} 
\begin{theorem}\label{prop:OO} 
	Let $\nu_1,\ldots,\nu_n$ be a basis for $\mathcal L$. Denote by $X_1,\ldots, X_n$ the corresponding basis of the Lie algebra $\gL$ of left invariant derivations of $\mathcal F$, and by $Y_1,\ldots, Y_n$ the corresponding basis of the Lie algebra $\gR$ of right invariant derivations of $\mathcal F$, that is, the ones such that $\epsilon \circ X_i = \epsilon \circ Y_i = \nu_i$ for every $i \in \{1,\ldots,n\}$. Denote the structure constants by $C_{ij}^k$ such that $[X_i,X_j] = \sum_k C_{ij}^k X_k$ for $i,j,k\in\{1,\ldots,n\}$. 
	
	Then there exist functions $$\OO^i_j, \bar\OO^i_j \in \mathcal F \cap \mathcal R, \quad i,j\in \{1,\ldots,n\} $$   such that $Y_j = \bar\OO^i_j X_i$ and $X_j = \OO^i_j  Y_i$ as elements of $\mathcal D$. Furthermore, these functions satisfy formulas 	
	\begin{equation}
	\epsilon (X_k(\OO^i_j)) = C^i_{kj},
	\end{equation}
	\begin{equation}\sum_{l,m} C_{lm}^k {\mathcal{O}}_i^l \mathcal{O}_j^m  = \sum_{r}  \mathcal{O}_r^k C_{ij}^r, \label{COOisOC}
	\end{equation}
	\begin{equation} 
	\sum_j \mathcal{O}_j^i \bar{\mathcal{O}}_k^j = \delta_k ^i = \sum_j \bar{\mathcal{O}}_j^i \mathcal{O}_k^j, \label{ObarOisI}
	\end{equation} 
	for $i,j,k\in \{1,\ldots,n\}$.
	Let $\mathcal H \subset \mathcal F$ be any Hopf algebra with the coalgebra structure and the antipode transposed from the group structure of $G$ such that it contains functions $\OO^i_j, \bar\OO^i_j, i,j\in\{1,\ldots,n\}$. Then the Hopf algebra $\mathcal H$ structure maps on them are
	\begin{equation} 
	\Delta(\mathcal{O}_j^i)= \sum_k \mathcal{O}_k^i \otimes \mathcal{O}_j^k, \quad \Delta(\bar{\mathcal{O}}_j^i)= \sum_k \bar{\mathcal{O}}_j^k \otimes \bar{\mathcal{O}}_k^i,  
	\end{equation} 
	\begin{equation}  
	%\vspace{-1em}
	\epsilon(\mathcal{O}_j^i)= \delta_j^i = \epsilon(\bar{\mathcal{O}}_j^i), 
	\end{equation} 
	\begin{equation} 
	S(\mathcal{O}_j^i) = \bar{\mathcal{O}}_j^i, \quad S(\bar{\mathcal{O}}_j^i) = \mathcal{O}_j^i,
	\end{equation} 
	for $i,j\in \{1,\ldots,n\}$.
	
\end{theorem}
\begin{proof}
	For case (I), denote by $\OO(g)$ the composition $({L_g}_* \circ {R_{g^{-1}}}_*)_g \colon T_g G \to T_g G$. Then $\OO(g) (Y_\nu) = X_\nu$ for every tangent vector  $\nu \in T_eG =\mathcal L$. The right invariant derivations $Y_1,\ldots,Y_n$ comprise a basis of the module $\mathcal D$ of all smooth vector fields over $\mathcal F$. In this basis, the matrix of the operator $\OO$ is precisely the matrix of the adjoint representation $\Ad$ in the basis $\nu_1,\ldots, \nu_n$. 	
	Indeed, we have for $\Ad_g \colon T_e G \to T_e G$, $\Ad_g = (L_{g*} \circ R_{g^{-1}*})_e$,
	\begin{align*}
	\Ad_g ((Y_j)_e) & = \Ad_g ((X_j)_e) = (L_{g*} \circ R_{g^{-1}*})_e ((X_j)_e) = (R_{g^{-1}*} \circ L_{g*})_e ((X_j)_e)  \\
	& = (R_{g^{-1}*})_g (( L_{g*})_e ((X_j)_e)) =  (R_{g^{-1}*})_g ((X_j)_g)   \\
	& = \sum_i (R_{g^{-1}*})_g (\OO^i_j(g) (Y_i)_g) = \sum_j\OO^i_j(g) (R_{g^{-1}*})_g ((Y_i)_g) \\
	& = \sum_i \OO^i_j(g) (Y_i)_e,
	\end{align*}
	and therefore, we conclude that 
	\begin{equation}\label{OAd}
	\OO^i_j(g) = [\Ad_g]^i_j.
	\end{equation}
	Note that we are only comparing matrix entries here; the matrix on the left is the matrix of a linear operator on $T_gG$ with respect to the basis  $(Y_1)_g, \ldots, (Y_n)_g$ and the matrix on the right is the matrix of a linear operator on $T_eG$ with respect to the basis $(Y_1)_e, \ldots, (Y_n)_e$, equal to the basis $(X_1)_e, \ldots, (X_n)_e$.

	Therefore, the functions $\OO^i_j$ that are the components of the matrix of $\OO$ in the basis $Y_1,\ldots,Y_n$ are representative functions, and $X_j = \sum_i \OO^i_j Y_i $ in $\mathcal D$. Denote  for every $g \in G$ the inverse of $\OO(g)$ by $\bar \OO(g)$. Then $Y_j = \sum_i \bar \OO^i_j X_i$.
	
	Since $\Ad_{gh} = \Ad_g \circ \Ad_h$ for every $g,h\in G$, we conclude that $$\OO^i_j(gh) = \sum_m \OO^i_m(g) \OO^m_j(h)$$ and therefore $\Delta(\OO^i_j) = \sum_m \OO^i_m\otimes \OO^m_j$. In a similar way, we conclude that $\Delta(\bar\OO^i_j) = \bar\OO^m_j \otimes \bar\OO^i_m$. Since the inverse of $\Ad_g$ is $\Ad_{g^{-1}}$, it is evident that $S(\OO^i_j) = \bar \OO^i_j $ and $S(\bar\OO^i_j) = \OO^i_j$. Since  $\epsilon(\OO^i_j) = \OO^i_j(e)$ and this is a component of matrix $\Ad_e$, we conclude that $\epsilon(\OO^i_j) = \delta^i_j$ and analogously that $\epsilon(\bar\OO^i_j) = \delta^i_j$. We calculate, for every $i,j,k \in \{1,\ldots,n\}$, 
	\begin{align*}
	\epsilon(X_k(\OO^i_j)) = &\lim_{t\to 0} \frac{d}{dt}\OO^i_j(\exp(t \nu_k)) = \lim_{t\to 0} \frac{d}{dt}[\Ad_{\exp(t \nu_k)}]^i_j \\
	& = \lim_{t\to 0} \,[\ad_{t\nu_k}]^i_j = C^i_{kj}.
	\end{align*}
	The property $C\OO\OO = \OO C$ in (\ref{COOisOC}) follows from the fact that $\Ad_g$ is a Lie algebra automorphism for every $g\in G$, with the identification $\gL \cong T_eG$, $X_i \mapsto \epsilon\circ X_i = (X_i)_e$, for $i\in\{1,\ldots,n\}$. Indeed, by writting $$\Ad_g ([X_i,X_j]_e) = [\Ad_g(X_i)_e , \Ad_g(X_j)_e], \quad \text{ for } i,j\in\{1,\ldots,n\}$$ in terms of the structure constants for $\gL$ and in terms of the matrix elements for $\Ad$ from formula (\ref{OAd}), we have that,  for every $g\in G$,
	\begin{align*}
	\Ad_g \Big(\sum_k C_{ij}^k (X_k)_e\Big) & = [\Ad_g(X_i)_e , \Ad_g(X_j)_e] \\
	\sum_k C_{ij}^k \OO_k^m(g) (X_m)_e & = \Big[\sum_l \OO^l_i(g) (X_l)_e, \sum_r \OO_j^r (g)(X_r)_e\Big] \\
	\sum_k C_{ij}^k \OO_k^m(g) (X_m)_e & = \sum_{l,r} \OO^l_i (g)\OO_j^r(g) C_{lr}^m (X_m)_e.
	\end{align*}The property $\OO \bar \OO = I = \bar\OO \OO$ in (\ref{ObarOisI}) follows from the fact that $\bar\OO(g)$ is a left and a right inverse to $\OO(g)$ for every $g \in G$. 
	
	For case (II), the unique functions exist by Propositions \ref{prop:isoleft} and \ref{prop:isoright} and they are representable because $\mathcal F \subset \mathcal R$. If we denote them by $\OO^ i_j$ and $\bar\OO^i_j$ such that $Y_j = \bar\OO^i_j X_i $ and $X_j = \OO^i_j Y_i$, we can prove the required equalities. Denote by $\OO$ and $\bar\OO$ matrices with components $\OO^i_j$ and $\bar\OO^i_j$ respectively. Since 
	\begin{align*}\sum_{m,i} \OO^i_m \bar\OO^m_j Y_i &= \sum_m \bar\OO^m_j X_m = Y_j = \delta^i_j Y_i, \\
	\sum_{m,i} \bar\OO^i_m \OO^m_j X_i & = \sum_m \OO^m_j Y_m = X_j = \delta^i_j X_i,
	\end{align*} for every $i \in \{ 1,\ldots, n \}$, we conclude that $\OO \bar \OO = I = \bar\OO \OO$, because the presentation of a derivation by a function and a basis of left (resp.\ right) invariant derivations is unique by Propositions \ref{prop:isoleft} and \ref{prop:isoright}. 
	Similarly, since in $\mathcal D$ we have 
	\begin{align*}
	[X_i, X_j] &= \sum_n C^n_{ij} X_n = \sum_{k,n} C^n_{ij} \OO^k_n Y_k,\\
	[X_i, X_j] &= [X_i, \sum_k \OO^k_j Y_k] \\
	&=   X_i \circ \sum_{k} \OO^k_j Y_k - \sum_k\OO^k_j  Y_k \circ X_i \\
	&= \sum_k X_i(\OO^k_j) Y_j + \sum_k \OO^k_j X_i \circ Y_j - \sum_k \OO^k_j Y_j \circ X_i \\
	&= \sum_k X_i(\OO^k_j) Y_j, \quad \text{ by Proposition } \ref{prop:comm},
	\end{align*} we conclude that $$X_i(\OO^k_j) = \sum_{n} C^n_{ij} \OO^k_n $$ and hence that 
	$$\epsilon(X_i(\OO^k_j)) = \langle \nu_i , \OO^k_j \rangle = C^k_{ij}.$$ 
	%	Furthermore, we know that 
	%	\begin{align*}
	%	[\sum_i \bar\OO^i_j X_i, X_m] 
	%	& = [Y_j, X_m] = 0,\\
	%	[\sum_i \bar\OO^i_j X_i, X_m] 
	%	&= [\sum_i \bar\OO^i_j X_i, X_m] = \sum_i \bar\OO^i_j X_i\circ  X_m - X_m \circ \sum_i \bar\OO^i_j X_i \\
	%	& = \sum_i \bar \OO^i_j X_i\circ X_m - \sum_i X_m(\bar\OO^i_j)X_i - \sum_i\bar\OO^i_j  X_m \circ X_i \\
	%	&= \sum_i \bar\OO^i_j [X_i,X_m] - \sum_i X_m(\bar\OO^i_j)X_i \\
	%	&= \sum_{i,k} \bar\OO^i_j C_{im}^k X_k - \sum_k X_m(\bar\OO^k_j)X_k ,
	%	\end{align*} we conclude that $$X_m(\bar\OO^k_j) = \sum_{i} C^k_{im} \bar\OO^i_j $$ and hence that 
	%	$$\epsilon(X_m(\bar\OO^k_j)) = C^k_{jm} = - C^k_{mj}.$$ 
	%In an analogous way, we can prove that $$Y_k(\OO^i_j) = \sum_n \tilde C^n_{kj} \OO^i_n.$$ ???
	It remains to prove $C\OO\OO = \OO C$. This follows from the identity 
	$$[X_i,X_j] = -\sum_{m,k}\OO_i^k \OO_j^m [Y_k,Y_m], \text{ for all } i,j \in \{1,\ldots,n\}$$ because this identity is equivalent to 
	$$C_{ij}^n \OO^k_n Y_k = C^k_{sm} \OO_i^s \OO_j^m Y_k, \text{ for all } i,j,k \in \{1,\ldots,n\},$$ 
	by using $[Y_i,Y_j] = -C_{ij}^k Y_k$ from Proposition \ref{minus}. 
	We have that 
	\begin{align*}
	[X_i,X_j] &=  [\sum_m \OO_i^m Y_m, \sum_k \OO^k_j Y_k] \\
	&= \sum_m \OO_i^m Y_m \circ \sum_k \OO^k_j Y_k -  \sum_k \OO^k_j Y_k \circ \sum_m \OO_i^m Y_m \\
	&= \sum_{m,k} \OO_i^m Y_m (\OO^k_j) Y_k  +  \sum_{m,k} \OO_i^m  \OO^k_j [Y_m, Y_k] -  \sum_{m,k} \OO^k_j Y_k (\OO_i^m) Y_m \\
	& = \sum_{k} X_i (\OO^k_j) Y_k +  \sum_{m,k} \OO_i^m  \OO^k_j [Y_m, Y_k] - \sum_{m} X_j (\OO_i^m) Y_m .
	\end{align*}	
	We have that
	\begin{align*}
	\sum_{k} X_i (\OO^k_j) Y_k  - \sum_{m} X_j (\OO_i^m) Y_m 
	&= \sum_{k,n} C^n_{ij}\OO_n^k Y_k  - \sum_{m,n} C^n_{ji}\OO^m_n Y_m \\
	& = \sum_{n} C_{ij}^n X_n   - \sum_{n} C^n_{ji}X_n \\
	& = [X_i,X_j] - [X_j,X_i] = 2[X_i,X_j],
	\end{align*}
	and therefore we conclude that 
	$$[X_i,X_j] + \sum_{m,k} \OO_i^m  \OO^k_j [Y_m, Y_k] = 0. $$ 
\end{proof}	

\begin{definition} Denote by $\OO^{min}(G)$ the smallest Hopf subalgebra of $\mathcal R \cap \mathcal F$ that contains functions $\OO^i_j, \bar\OO^i_j ,i,j\in \{1,\ldots,n\} $ from Theorem \ref{prop:OO}, in both cases (I) and~(II).
\end{definition}

\subsection{Left coaction}

\begin{theorem}\label{prop:lambdaO} The map $\gL \to \mathcal L \to \gR \hookrightarrow \mathcal F \otimes \gL$, that sends the left invariant derivation to the corresponding right invariant one written as an element of $\mathcal F \otimes  \gL$, extends to a unique antimultiplicative unital linear map $\lambda \colon \UgL \to  \mathcal F \sharp_{\str} \UgL$.
	
Furthermore, for every Hopf algebra $\mathcal H$ such that $\OO^{min}(G) \subset \mathcal H  \subset \mathcal F \cap\mathcal R$, the corestriction $$\lambda \colon \UgL \to \mathcal H \sharp_{\str} \UgL$$ exists and it is a left coaction. The elements of $\Img \lambda$ commute with the elements of $1\sharp \UgL$ in $\mathcal H \sharp \UgL$.
	
\end{theorem}
\begin{proof}

		(i) We prove that such $\lambda$ exists. We first define an auxiliary map $\tilde\lambda$ as a linear map $\gL \to \mathcal H \sharp U(\gL)$ such that $\tilde\lambda (X_j) = \sum_i \bar{\mathcal{O}}_j^i \sharp X_i$ for $j \in \{1,\ldots,n\}$ and expand it to $\tilde \lambda \colon T(\gL) \to \mathcal H \sharp U(\gL)$ by antimultiplicativity, and then we check that $\tilde\lambda([X,Z] - X\otimes X' + X'\otimes X) = 0$ for every $X,X' \in \gL$. We have 
		\begin{align*} \tilde\lambda(X_j) \cdot \tilde\lambda(X_i) &= (\sum_k \bar{\mathcal O}_j^k\sharp X_k) \cdot( \sum_m \bar{\mathcal O}_i^m \sharp X_m) = \sum_{k,m,p} \bar{\mathcal O}_j^k \bar{\mathcal O}^p_i (X_k \btl \bar{\mathcal O}_p^m) X_m \\
		& = \sum_{k,m,p} \bar{\mathcal O}_j^k \bar{\mathcal O}^p_i \sharp (\delta_p^m X_k  - C_{kp}^m) X_m   \\
		&= \sum_{k,m} \bar{\mathcal O}_j^k \bar{\mathcal O}^m_i \sharp X_k X_m - \sum_{k,m,p} \bar{\mathcal O}_j^k \bar{\mathcal O}^p_i  C_{kp}^m\sharp X_m  
		\end{align*}
		and, analogously, we have
		\begin{align*}
		\tilde\lambda(X_i) \cdot \tilde\lambda(X_j) &= \sum_{k,m} \bar{\mathcal O}^m_i \bar{\mathcal O}^k_j \sharp X_m X_k  - \sum_{k,m,p} \bar{\mathcal O}^m_i \bar{\mathcal O}^p_j  C_{mp}^k\sharp X_k.
		\end{align*}  
		On the other hand, we have
		$$\tilde\lambda([X_i,X_j]) = \tilde\lambda(\sum_p C^p_{ij} X_p) = \sum_{p,m} C_{ij}^p \bar{\mathcal O}_{p}^m \sharp X_m.$$
		Equality $\tilde\lambda(X_j)\tilde\lambda(X_i) - \tilde\lambda(X_i)\tilde\lambda(X_j) = \tilde\lambda([X_i,X_j])$ now follows from commutativity of $\mathcal H$ and identities $$\sum_{l,m} C_{lm}^k \bar{\mathcal{O}}_i^l \bar{\mathcal{O}}_j^m = \sum_{r} \bar{\mathcal{O}}_r^k C_{ij}^r , \qquad  i,j, k \in \{1,\ldots,n\},$$ 
		which are easily deduced from   
		$$\sum_{l,m} C_{lm}^k{\mathcal{O}}_i^l \mathcal{O}_j^m = \sum_{r}  \mathcal{O}_r^k C_{ij}^r, \qquad  i,j, k\in \{1,\ldots,n\}.$$ Therefore, by quotienting the domain of $\tilde \lambda$ with the ideal generated by all elements of the form $[X,X'] - X\otimes X' + X'\otimes X$, where $X\in \gL$, we induce a well defined map $\lambda \colon U(\gL) \to \mathcal H \sharp U(\gL)$.
		Additionally, we note that clearly 
		\begin{equation}\label{eq:antimult} 
		\lambda(DP) = \lambda(P)\lambda(D),\quad \text{ for all } D,P \in U(\gL).
		\end{equation}
		(ii) We prove that the elements of $\operatorname{Im}\lambda $ commute with the elements of $1\sharp \UgL$ in $\mathcal H\sharp \UgL$.
		First we check this for generators. For any $j,k \in \{1,\ldots,n\}$ we have
		\begin{align*}
		X_k \cdot \lambda(X_j) &=  X_k \cdot \sum_i  \bar{\mathcal O}_j^i \sharp X_i = \sum_{i,m} \bar\OO^m_j \sharp(X_k \btl \bar{\mathcal O}_m^i ) X_i \\
		&=		\sum_{i,m} \bar\OO^m_j \sharp(\delta_m^i X_k + C_{mk}^i)X_i = \sum_{m} \bar\OO^m_j \sharp(X_k X_m + [X_m,X_k]) = \\
		&= 	\sum_{m} \bar\OO^m_j \sharp X_m X_k = \lambda(X_j) \cdot X_k.
		\end{align*} Then by using antimultiplicativity (\ref{eq:antimult}) of $\lambda$, it is easy to prove the claim inductively for all elements of $\UgL$,
		\begin{equation}\label{eq:imacom}
		P \cdot \lambda(D) = \lambda(D) \cdot P, \quad \text{ for all } D,P \in \UgL.
		\end{equation}

		(iii) We prove that $\lambda$ is a coaction. First, we prove that if the coassociativity identity is true for  $D,P \in U(\gL)$, then it is true for the product $DP \in U(\gL)$. We use these properties: $\lambda$ is an antihomomorphism of algebras as proven in (i), $\btl$ is a Hopf action, properties (\ref{eq:imacom}) proven in (ii), and compatibility of comultiplication and multiplication of $\mathcal H$. We thus have that
		\begin{align}\lambda(DP) &= \lambda(P) \lambda(D) = \sum P_{[-1]}\sharp P_{[0]} \cdot \lambda(D) = \nonumber \\
		&= \sum P_{[-1]} \cdot \lambda(D) \cdot P_{[0]} \quad \text{ by (\ref{eq:imacom})} \nonumber\\
		&= \sum P_{[-1]} D_{[-1]} \sharp D_{[0]} P_{[0]}  \label{eq:coabc}
		\end{align}
		from which it follows, assuming the coassociativity identity holds for $D$ and $P$,
		\begin{align*}
		((\id \otimes \lambda)\circ \lambda)(PD) 
		& = \sum D_{[-1]} P_{[-1]} \otimes \lambda( P_{[0]} D_{[0]})  \\
		&=  \sum D_{[-1]} P_{[-1]} \otimes D_{[0][-1]} P_{[0][-1]} \otimes P_{[0][0]} D_{[0][0]} \\
		& = \sum D_{[-1](1)} P_{[-1](1)} \otimes D_{[-1](2)} P_{[-1](2)} \otimes P_{[0]} D_{[0]} \\
		& = \sum (D_{[-1]} P_{[-1]})_{(1)} \otimes (D_{[-1]}P_{[-1]})_{(2)}  \otimes P_{[0]} D_{[0]}
		\end{align*}
		and, on the other hand,  
		\begin{align*}
		((\Delta \otimes \id) \circ \lambda)(PD) & = (\Delta \otimes \id)(\sum D_{[-1]} P_{[-1]} \sharp P_{[0]} D_{[0]}) \\
		& = \sum (D_{[-1]}P_{[-1]})_{(1)} \otimes (D_{[-1]}P_{[-1]})_{(2)}  \otimes P_{[0]} D_{[0]}.
		\end{align*}
		Now we see that it is sufficient to check the claim on generators, and this is trivial:
		$$\sum_{k,i} \bar{\mathcal O}_j^k \otimes \bar{\mathcal O}_k^i \sharp X_i = \sum_{i,m} \bar{\mathcal O}_j^i \otimes \bar{\mathcal O}_i^m \sharp X_m.$$
		This proves the coassociativity of $\lambda$. The counitality of $\lambda$ is checked similarly, first on generators, for every $j\in \{1,\ldots,n\}$, $((\epsilon \otimes \id) \circ \lambda)(X_j) = \epsilon(\bar\OO^i_j) X_i = X_j,$ and then easily proven inductively by using formula (\ref{eq:coabc}).

\end{proof}

\begin{theorem}\label{prop:lambdaOR} The map $\gR \to \mathcal L \to \gL \hookrightarrow \mathcal F \otimes \gR$, that sends the right invariant derivation to the corresponding left invariant one written as an element of $\mathcal F \otimes  \gR$, extends to a unique antimultiplicative unital linear map $\lambda \colon \UgR \to  \mathcal F \sharp_{\str} \UgR$.
	
Furthermore, for every Hopf algebra $\mathcal H$ such that $\OO^{min}(G) \subset \mathcal H  \subset \mathcal F$, the corestriction $$\lambda \colon \UgR \to \mathcal H \sharp_{\str} \UgR$$ exists and it is a left coaction. The elements of $\Img \lambda$ commute with the elements of $1\sharp \UgL$ in $\mathcal H \sharp \UgL$.
\end{theorem}
\begin{proof}
	Analogous to the proof of Theorem \ref{prop:lambdaO}.
\end{proof}

\subsection{Right action that produces the same smash product algebra}

\begin{proposition}\label{prop:rightact} Maps  
	$$\btl \colon \UgL \otimes \mathcal H \to \UgL, \quad \btl \colon \UgR \otimes \mathcal H^\co \to \UgR$$ defined by the same formula 
	$$D \btl f = \langle D_{(1)}, f \rangle D_{(2)}, \quad \text{ for } f \in \mathcal H \text{ and } D \in \UgL \text{ or } \UgR,$$ are right Hopf actions and they induce smash product algebras isomorphic to the ones induced by the left actions in Proposition \ref{prop:smash},
	$$\mathcal H \sharp_{\str} \UgL \cong \mathcal H \sharp_{\btl} \UgL, \quad \mathcal H \sharp_{\str} \UgR \cong \mathcal H^\co \sharp_{\btl} \UgR.$$

\end{proposition}
\begin{proof} See Proposition \ref{prop:leftright}.
\end{proof}
\subsection{Yetter--Drinfeld property and the main theorem}

\begin{theorem}\label{prop:YDO} Let a Hopf algebra $\mathcal H$ be a subalgebra of $\mathcal F$ with the coalgebra structure and the antipode transposed from the group structure on $G$ and assume that $\OO^{min}(G) \subset \mathcal H$. 
	
	Then $(\UgL, \btl, \lambda)$ is a  braided commutative Yetter--Drinfeld module algebra over~$\mathcal H$. Similarly, $(\UgR, \btl, \lambda)$  is a  braided commutative Yetter--Drinfeld module algebra over~$\mathcal H^\co$. Here symbols $\btl$ and $\lambda$ denote the appropriate action and coaction in each case, from Proposition \ref{prop:rightact} and Theorems \ref{prop:lambdaO} and \ref{prop:lambdaOR}.
	
\end{theorem}
\begin{proof} (i) We prove the Yetter--Drinfeld property for $(\UgL, \btl, \lambda)$,
	$$\sum f_{(2)} \cdot \lambda(D\btl f_{(1)}) = \lambda(D) \cdot f, \quad \text{ for all } D\in \UgL \text{ and } f\in \mathcal H.$$ 
	It is sufficient to prove this for generators  $X_k, 1_{\UgL} \in \UgL$, because $\btl$ is a Hopf action and $\lambda$ is antimultiplicative. Indeed, if we assume that the equality holds for some $D$ and $P$ in $\UgL$ and every $f\in\mathcal H$, then it also holds for the product $DP$ and all $f\in\mathcal H$, by
	\begin{align*}
	\sum f_{(2)} \lambda((DP) \btl f_{(1)}) & = \sum f_{(3)} \lambda((D \btl f_{(1)})(P \btl f_{(2)})) \\
	& = \sum f_{(3)} \lambda(P \btl f_{(2)})\lambda(D \btl f_{(1)}) \\
	& =\sum\lambda(P)\cdot f_{(2)} \lambda(D \btl f_{(1)}) \\
	&= \lambda(P) \lambda(D) \cdot f \\
	&= \lambda(DP) \cdot f.
	\end{align*}
	The check for $1_{\UgL}$ is trivial. Let us prove the claim now for generators:	
	$$f_{(2)} \lambda(X_k \btl f_{(1)}) = \lambda(X_k) \cdot f, \quad \text{ for every } k \in \{1,\ldots,n\} \text{ and every } f \in \mathcal H.$$ By the definition of the right action of $\mathcal H$ on $\UgL$, this is $$f_{(2)}\langle (X_k)_{(1)}, f_{(1)} \rangle \lambda((X_k)_{(2)}) = \lambda(X_k) \cdot f,$$ and since $\Delta(X_k) = X_k \otimes 1 + 1 \otimes X_k$ and the pairing is Hopf, we have that this is $$f_{(2)}\langle X_k, f_{(1)} \rangle \lambda(1) + f_{(2)}\langle 1, f_{(1)} \rangle \lambda(X_k) = \lambda(X_k) \cdot f,$$ that is $$f_{(2)}\langle X_k, f_{(1)} \rangle + f \lambda(X_k) = \lambda(X_k) \cdot f.$$ 
	Since $\langle X_k, f \rangle = \langle Y_k, f \rangle$ for all $f \in \mathcal H$, this is equivalent to $$f_{(2)}\langle Y_k, f_{(1)} \rangle + f \lambda(X_k) = \lambda(X_k) \cdot f,$$ that is $$Y_k \str f = Y_k \cdot f - f \cdot Y_k,$$ which holds because this left action is a differentiation by a right invariant derivation.

	(ii) The comodule algebra property follows from the antimultiplicativity of $\lambda$ and the commuting of elements of $\Img \lambda$ with elements of $1\sharp \UgL$, proven in Theorem~\ref{prop:lambdaO}.  
	
	(iii) The braided commutativity property is: 
	\begin{equation}\label{eq:lambdact} P \btl \lambda(D) = DP, \quad \text{ for all } D,P \in U(\gL).
	\end{equation} First we check this on generators. For any $j,k \in \{1,\ldots,n\}$ we have
	$$X_k \btl \sum_i \bar{\mathcal O}_j^i \sharp X_i = \sum_i (\delta_j^i X_k X_i + C_{jk}^i X_i) = X_k X_j + [X_j,X_k] = X_j X_k.$$ 
	Next, we use induction on the length of the word acted on by $\lambda(X_j)$ on the right, for every $X_j, j\in\{1,\ldots,n \}$. The step of induction is: 
	$$(PD) \btl \sum_i \bar{\mathcal O}_j^i \sharp X_i = \sum_{i,m} (P  \btl  \bar{\mathcal O}_j^m)(D \btl \bar{\mathcal O}_m^i )X_i = \sum_{i,m} (P  \btl  \bar{\mathcal O}_j^m)X_m D =  X_j PD.$$ 
	At last, we use induction on the length of the word on the right. The step of induction is:
	$$D\btl \lambda(PT) = (D \btl \lambda(T)) \btl \lambda(P) = (TD) \btl \lambda(P) = PTD.$$
	
	The proof of (i), (ii) and (iii) for $(\UgR, \btl,\lambda)$ is analogous.
\end{proof}

\begin{corollary}
	Let a Hopf algebra $\mathcal H$ be a subalgebra of $\mathcal F$ with the coalgebra structure and the antipode transposed from the group structure on $G$ and assume that $\OO^{min}(G) \subset \mathcal H$. 
	
	Then $\mathcal H \sharp \UgL \cong \UgR \sharp \mathcal H$ is a scalar extension Hopf algebroid over base algebras $\UgR, \UgL$ and the isomorphic algebra $\mathcal H^\co \sharp \UgR \cong \UgL \sharp \mathcal H^\co$ is a scalar extension Hopf algebroid over base algebras $\UgL, \UgR$
	
\end{corollary}
\begin{proof}
	See Remark 4.1 and Theorem 4.2 in \cite{stojic}.
\end{proof}

\section{Yetter--Drinfeld module algebra over a Hopf algebra of functionals}

Here we consider case (III) described in Section \ref{sec:cases}.  Let now $\ggf$ be any finite-dimesional Lie algebra over any field $\genfd$. Denote by $\Ug^*$ the dual of the universal enveloping algebra $\Ug$. The Hopf algebra structure on the finite dual $\Ug^\circ \subset \Ug^*$ is such that the restriction $\Ug \otimes \Ug^\circ \to \genfd$ of the canonical pairing $\Ug \otimes \Ug^* \to \genfd$ is a Hopf pairing.

\subsection{Definition of Hopf algebra $\Ug^{min}$}

We now choose a basis of $\ggf$ and define matrices $\UU \in M_n(\Ug^\circ)$ and $\bar{\UU} \in M_n(\Ug^\circ)$ such that, if there were an affine algebraic group whose Lie algebra is $\ggf$, functionals that are components of $\UU$ would pair with $\Ug$ precisely as functions that are components of~$\OO$, that is, the smallest Hopf algebra $\Ug^{min}$ generated by these functionals would be the image of $\OO^{min}(G)$ in $\Ug^*$ under mapping $\OO(G) \to \Ug^*$ defined by $f \mapsto \langle - , f \rangle$ for $f\in \OO(G)$.
$$
\begin{tikzcd}
\Ug^* & \ar{l}{} \OO(G)  & \\
\Ug^{min} \ar[hook]{u}  &  \ar[two heads]{l} \OO^{min}(G) \ar[hook]{u}  & \text{ for }\ggf = \gL & &
\end{tikzcd}
$$
All maps in the diagram are Hopf algebra morphisms that commute with the pairings with $\Ug$, provided that we assume the completed version of a Hopf algebra structure on $\Ug$. Similarly, for the right invariant case, we have the diagram
$$
\begin{tikzcd}
\Ug^* & \ar{l}{} \OO(G)^\co  & \\
\Ug^{min} \ar[hook]{u}  &  \ar[two heads]{l} \OO^{min}(G)^\co \ar[hook]{u}  & \text{ for }\ggf = \gR. & &
\end{tikzcd}
$$
We do this by requiring that $\langle X,\UU\rangle$ is the matrix presentation of $\ad_X$ with respect to the chosen basis and that $\langle X,\bar\UU\rangle$ is the matrix presentation of $-\ad_X$ with respect to the chosen basis, for every $X \in \ggf$, and such that certain equations regarding the comultiplication of their components hold. 

\begin{theorem}\label{prop:U}
	Let $X_1,\ldots, X_n$ be a basis for a Lie algebra $\ggf$ over a field $\genfd$. Let the structure constants $C^i_{jk}$ be defined from $[X_i,X_j] = \sum_k C^i_{jk} X_k$, for $i,j\in \{1,\ldots,n\}$.
	
	Then there exist unique functionals $$\UU^i_j, \bar\UU^i_j \in \Ug^\circ, \quad i,j\in \{1,\ldots,n\} $$   such that, for all $i,j,k \in \{1,\ldots,n\}$, 	
	$$\langle X_k, \UU^i_j \rangle = C^i_{kj},  $$
	$$\Delta(\mathcal{U}_j^i)= \sum_k \mathcal{U}_k^i \otimes \mathcal{U}_j^k, \quad \Delta(\bar{\mathcal{U}}_j^i)= \sum_k \bar{\mathcal{U}}_j^k \otimes \bar{\mathcal{U}}_k^i,   $$ 
	\vspace{-1em}
	$$ \epsilon(\mathcal{U}_j^i)= \delta_j^i = \epsilon(\bar{\mathcal{U}}_j^i), $$
	$$S(\mathcal{U}_j^i) = \bar{\mathcal{U}}_j^i, \quad S(\bar{\mathcal{U}}_j^i) = \mathcal{U}_j^i.$$
	These functionals satisfy the following identities, for all $i,j,k\in \{1,\ldots,n\}$,
	\begin{align}
	&\sum_{l,m} C_{lm}^k {\mathcal{U}}_i^l \mathcal{U}_j^m = \sum_{r}  \mathcal{U}_r^k C_{ij}^r, \label{ccu}\\
	&\sum_j \mathcal{U}_j^i \bar{\mathcal{U}}_k^j = \delta_k ^i = \sum_j \bar{\mathcal{U}}_j^i \mathcal{U}_k^j. \label{cui}
	\end{align}
\end{theorem} 
\begin{proof} 
	
	\emph{Uniqueness.}  If matrices $\UU,\bar\UU \in M_n(\Ug^\circ)$ with these properties exists, then they are unique since for all $i,j\in \{1,\ldots,n\}$ we have that 
	$$\langle X_{i_1} \cdots X_{i_m}, \UU^i_j  \rangle = \langle  X_{i_1} \otimes \cdots \otimes X_{i_m}, \Delta_{m}(\UU^i_j) \rangle, \quad \text{ for all } X_{i_1},\ldots,X_{i_n} \in \ggf $$ and comultiplication of $\UU^i_j$ is written in terms of tensor products of these functionals. Next, we have for all $i,j\in \{1,\ldots,n\}$ 
	$$\langle  X_{i_1} \cdots X_{i_m}, \bar\UU^i_j \rangle = \langle  S(X_{i_1} \cdots X_{i_m} ), \UU^i_j\rangle, \quad \text{ for all } X_{i_1},\ldots,X_{i_n} \in \ggf .$$

	\emph{Existence and formulas.} Denote by $\mathcal{C}_k = [C_{kj}^i]^i_j$ the matrix of $\ad_{X_k}$ in this basis. 
	We define $\UU$ and $\bar\UU$ first as elements of $M_n(T(\ggf)^*)$ by
	$$\langle 1, \UU \rangle = I = \langle 1, \bar\UU \rangle$$ 
	$$\langle X_{i_1} \cdots X_{i_m}, \UU \rangle = \mathcal{C}_{i_1} \cdots \mathcal{C}_{i_m}, \quad \langle X_{i_1} \cdots X_{i_m}, \bar\UU \rangle = (-1)^{m} \mathcal{C}_{i_m} \cdots \mathcal{C}_{i_1}$$ 
	and then check that they are well defined as elements of $M_n(\Ug^*)$. We have that
	\begin{align*}
	\langle [X_i,X_j], \UU^k_m \rangle &= \langle \sum_s C_{ij}^s X_s,\UU^k_m \rangle = \sum_s C_{ij}^s C^k_{sm}.
	\end{align*}
	Next, we have that
	\begin{align}
	\sum_k \langle X_iX_j - X_j X_i, \mathcal U^k_m \rangle X_k 
	&= \sum_{k,s} C_{is}^k C_{jm}^s X_k - C^k_{js}  C^s_{im} X_k \nonumber\\
	& = \sum_{s} [X_i,X_s] C_{jm}^s- \sum_s [X_j,X_s]   C^s_{im} \nonumber \\
	& = [X_i,[X_j,X_m]] - [X_j,[X_i,X_m]] \nonumber \\ 
	&= [[X_i,X_j], X_m] = \sum_{s,k} C_{ij}^s C_{sm}^k X_k \label{cc}
	\end{align}
	from which we conclude that $$\langle X_iX_j - X_j X_i, \mathcal U^k_m \rangle = \sum_{s} C_{ij}^s C_{sm}^k.$$
	The claim for generators $\bar\UU^i_j$ is proven analogously. Therefore, they are well defined as elements of $\Ug^*$. It is easy to verify that the Hopf algebra structure maps indeed are as should be. Therefore, $\UU, \bar\UU \in M_n(\Ug^\circ)$. 
	
	\emph{Identities}. The identity (\ref{ccu}) is proven as follows. For every $i,j,k,s \in \{1,\ldots,n\}$
	\begin{align*}
	\langle  X_s, \sum_{l,m} C_{lm}^k {\mathcal{U}}_i^l \mathcal{U}_j^m \rangle 
		&= \sum_{l,m} C_{lm}^k  \langle X_s, {\mathcal{U}}_i^l  \rangle \delta_j^m + \sum_{l,m} C_{lm}^k \delta_i^l \langle    X_s , \mathcal{U}_j^m \rangle \\
		&= \sum_{l} C_{lj}^k  C^l_{si} + \sum_{m} C_{im}^k  C_{sj}^m, \\
	\langle \sum_{r}  \mathcal{U}_r^k C_{ij}^r ,X_s \rangle 
		&= \sum_{r}  C^k_{rs}C^r_{ij}.
		\end{align*}
	The equality is now proven similarly as in calculation (\ref{cc}). 

\end{proof}	

\begin{definition} Denote by  $\Ug^{min}$ the smallest subalgebra of $\Ug^*$ generated by the components of the matrices $\UU$ and $\bar\UU$ from Theorem \ref{prop:U}. Clearly $\Ug^{min}$ is a Hopf subalgebra of $\Ug^\circ$.
\end{definition}

\subsection{Left coaction}

\begin{proposition}\label{prop:lambdaU} Let $\ggf$ be a finite-dimensional Lie algebra over a field $\genfd$. Let $\mathcal H$ be a Hopf $\genfd$-algebra in a Hopf pairing with $\Ug$ such that $\Ug^{min} \subset \mathcal H \subset \Ug^\circ$. 
	Then the non-degenerate pairing $\Ug \otimes\mathcal  H \to \genfd$ induces a unique right Hopf action $\btl \colon \Ug \otimes\mathcal  H \to \Ug$ such that $\langle D\btl f, k \rangle = \langle D, fk\rangle$ for all $D \in \Ug$ and all $f,k \in\mathcal  H$ and this action induces a smash product algebra structure $\mathcal H \sharp \Ug$. 
	
	Let $X_1,\ldots, X_n$ be a basis of $\ggf$. Let $\UU$ and $\bar\UU \in M_n(\Ug^\circ)$ be as in Theorem \ref{prop:U}, with regard to this basis. 
	
	Then there exists a unique linear antimultiplicative unital map $$\lambda \colon \Ug \to \mathcal H \sharp \Ug \subset \Ug^\circ \sharp \Ug  $$ such that $\lambda(X_j) = \sum \bar\UU^i_j \sharp X_i $ for every $j \in \{1,\ldots,n\}$. The elements of  $\Img \lambda$ commute with the elements of $1\sharp \Ug$ in $\mathcal H\sharp \Ug$.
\end{proposition}
\begin{proof} The proof is essentially the same as the proof of Theorem \ref{prop:YDO}. 

\end{proof}
\begin{proposition}\label{prop:lambdaUR} Let $\ggf$ be a finite-dimensional Lie algebra over a field $\genfd$. Let $\mathcal H$ be a Hopf $\genfd$-algebra in a Hopf pairing with $\Ug$ such that $\Ug^{min} \subset \mathcal H \subset \Ug^\circ$. 
	Then the non-degenerate pairing $\Ug \otimes\mathcal  H \to \genfd$ induces a unique right Hopf action $\btl \colon \Ug \otimes\mathcal  H \to \Ug$ such that $\langle D\btl f, k \rangle = \langle D, fk\rangle$ for all $D \in \Ug$ and all $f,k \in\mathcal  H$ and this action induces a smash product algebra structure $\mathcal H \sharp \Ug$. 
	
	Let $Y_1,\ldots, Y_n$ be a basis of $\ggf$. Let $\UU$ and $\bar\UU \in M_n(\Ug^\circ)$ be as in Theorem \ref{prop:U}, with regard to this basis. 
	
	Then there exists a unique linear antimultiplicative unital map $$\lambda \colon \Ug \to \mathcal H \sharp \Ug \subset \Ug^\circ \sharp \Ug  $$ such that $\lambda(Y_j) = \sum \UU^i_j \sharp Y_i $ for every $j \in \{1,\ldots,n\}$. The elements of  $\Img \lambda$ commute with the elements of $1\sharp \Ug$ in $\mathcal H\sharp \Ug$.
\end{proposition}
\begin{proof}
	The proof is analogous to the proof of Theorem \ref{prop:lambdaU}.
\end{proof}
\subsection{Yetter--Drinfeld property and the main theorem}
\begin{theorem} Let $\ggf$ be a finite-dimesional Lie algebra over any field $\genfd$. Let $\mathcal H$ be any Hopf algebra such that $\Ug^{min} \subset \mathcal H \subset \Ug^\circ$. 
	Then the non-degenerate pairing $\Ug \otimes\mathcal  H \to \genfd$ induces a unique right Hopf action $\btl \colon \Ug \otimes\mathcal  H \to \Ug$ such that $\langle D\btl f, k \rangle = \langle D, fk\rangle$ for all $D \in \Ug$ and all $f,k \in\mathcal  H$ and this action induces a smash product algebra structure $\mathcal H \sharp \Ug$. 
	
	Let $\lambda \colon \Ug \to\mathcal  H \sharp \Ug$ be as defined in Proposition \ref{prop:lambdaU}. Then $(\Ug,\btl, \lambda)$ is a braided commutative Yetter--Drinfeld module algebra over $\mathcal H$.
\end{theorem}
\begin{proof} (i) We prove the Yetter--Drinfeld property:
	$$\sum f_{(2)} \cdot \lambda(X\btl f_{(1)}) = \lambda(X) \cdot f, \quad \text{ for all } X\in \UgL \text{ and } f\in\mathcal  H.$$
	
	It is sufficient to prove this for generators  $X_k, 1_{\Ug} \in \Ug$ because we have that: $\btl$ is a Hopf action and $\lambda$ is an antihomomorphism of algebras. Indeed, if we assume that the equality holds for $X$ and $X'$ in $\Ug$, then it also holds for $XX'$:
	\begin{align*}
	\sum f_{(2)} \lambda((XX') \btl f_{(1)}) & = \sum f_{(3)} \lambda((X \btl f_{(1)})(X' \btl f_{(2)})) \\
	& = \sum f_{(3)} \lambda(X' \btl f_{(2)})\lambda(X \btl f_{(1)}) \\
	& =\sum\lambda(X')\cdot f_{(2)} \lambda(X \btl f_{(1)}) \\
	&= \lambda(X') \lambda(X) \cdot f \\
	&= \lambda(XX') \cdot f.
	\end{align*}
	The check for $1_{\Ug}$ is trivial. Let us prove the claim now for generators:	
	$$f_{(2)} \lambda(X_k \btl f_{(1)}) = \lambda(X_k) \cdot f, \quad \text{ for every } k \in \{1,\ldots,n\} \text{ and every } f \in H.$$ 
	Extend the Hopf action $\str \colon \Ug \otimes\mathcal  H \to\mathcal  H$ defined by $D \str f = f_{(1)} \langle D,f_{(2)}\rangle$, to an action of the whole smash product algebra $\str \colon\mathcal  H \sharp \Ug \otimes\mathcal  H \to\mathcal  H $ by formula $f \sharp D \str g := f \cdot (D\str g)$ for $f,g \in \mathcal  H$ and $D \in \Ug$. The action $\btl \colon \Ug \otimes\mathcal  H \to \Ug$ also extends to an action of the whole smash product algebra $\btl \colon \Ug \otimes\mathcal  H\sharp \Ug \to \Ug$ by formula $X \btl f \sharp Z = (X \btl f) Z$. It is proven in Proposition \ref{prop:lambdaU} that the elements of $\Img \lambda$ commute with the elements of $1\sharp \Ug$ in $\mathcal H \sharp \Ug$. The pairing of $\Ug$ with $\mathcal H$ has the property $\langle D , f \rangle = \langle 1, D \str f \rangle$, and we extend it to the pairing $\mathcal H \sharp \Ug \otimes\mathcal  H \to \genfd$ by formula $\langle f\sharp D , g \rangle = \langle 1, (f\sharp D)\str g\rangle$, for for $f,g \in \mathcal  H$ and $D \in \Ug$. 
	Now we check the Yetter--Drinfeld property equation for $X\in \ggf$ and $f \in\mathcal  H$,
	$$\langle X, f_{(1)}\rangle f_{(2)} = \lambda(X) \cdot f - f \cdot \lambda(X).$$ 
	First, we check that $$\lambda(X) \str f = \lambda(X) \cdot f - f \cdot \lambda(X), \quad \text{ for all }X \in \ggf\text{ and }	f\in\mathcal  H,$$ 
	which is straightforward and uses the  commutativity of $\mathcal H$. 
	 Next, we prove that $$\langle X, f_{(1)}\rangle f_{(2)} = \lambda(X) \str f, \quad \text{ for all }X \in \ggf \text{ and } 
	f\in\mathcal  H. $$ Since the pairing is non-degenerate in the second variable, this is equivalent to  
	$\langle D, \langle X, f_{(1)} \rangle f_{(2)} \rangle = \langle D ,\lambda(X) \str f \rangle$ being true for all $X\in\ggf$, $ D\in\Ug$ and $f\in\mathcal  H$. The equation is $\langle X, f_{(1)}\rangle \langle D, f_{(2)} \rangle = \langle D \lambda(X), f \rangle $, that is $\langle XD, f \rangle = \langle D \lambda(X), f \rangle$. Since the elements of the image of $\lambda$ commute with the elements of $\Ug$, this is then $\langle XD, f \rangle = \langle \lambda(X) D, f \rangle$. For $X \in \ggf$ it is easy to check that $\langle XD, f \rangle = \langle \lambda(X)D, f \rangle$, by using $\epsilon(\UU) = I$.  
	
	(ii) The comodule algebra property and (iii) the braided commutativity property are proven in the same way as those properties in the proof of Theorem \ref{prop:YDO}.
\end{proof}

\begin{corollary}
	Let $\mathcal H$ be a subalgebra of $\UgL^*$ that is in a Hopf pairing with $\UgL$ and assume that $\UgL^{min} \subset \mathcal H$. 
	
	Then $\mathcal H \sharp \UgL \cong \UgR \sharp \mathcal H$ is a scalar extension Hopf algebroid over base algebras $\UgR, \UgL$ and the isomorphic algebra $\mathcal H^\co \sharp \UgR \cong \UgL \sharp \mathcal H^\co$ is a scalar extension Hopf algebroid over base algebras $\UgL, \UgR$.
	
\end{corollary}
\begin{proof}
	See Remark 4.1 and Theorem 4.2 in \cite{stojic}.
\end{proof}

\section{Theorem with all intermediate algebras} 
\begin{theorem} \label{prop:general}
	Let $\ggf$ be a finite-dimensional Lie algebra over a field $\genfd$, with structure constants $C^k_{ij}$ determined from $$[X_i,X_j] = C^k_{ij} X_k, \quad \text{ for } i,j\in\{1,\ldots,n\}.$$ Let $\mathcal H$ be a Hopf $\genfd$-algebra in a  Hopf pairing $\langle\,,\rangle \colon \Ug \otimes \mathcal H \to \genfd$ with $\Ug$. 
	
	If there exist some elements $\mathcal{O}^i_j, \bar{\mathcal{O}}^i_j, i,j\in \{1,\ldots,n\}$ of $\mathcal H$  such that %they generate all elements of $H$ as an algebra and such that	
	$$\langle X_k, \OO^i_j \rangle = C^i_{kj}, $$
	$$\sum_{l,m} C_{lm}^k {\mathcal{O}}_i^l \mathcal{O}_j^m = \sum_{r}  \mathcal{O}_r^k C_{ij}^r, $$ 
	$$\sum_j \mathcal{O}_j^i \bar{\mathcal{O}}_k^j = \delta_k ^i = \sum_j \bar{\mathcal{O}}_j^i \mathcal{O}_k^j, \quad \text{ for all } i,j,k\in \{1,\ldots,n\}$$
	and such that the coalgebra structure maps and the antipode of $\mathcal H$ satisfy 
	$$\Delta(\mathcal{O}_j^i)= \sum_k \mathcal{O}_k^i \otimes \mathcal{O}_j^k, \quad \Delta(\bar{\mathcal{O}}_j^i)= \sum_k \bar{\mathcal{O}}_j^k \otimes \bar{\mathcal{O}}_k^i,   $$ 
	\vspace{-1em}
	$$ \epsilon(\mathcal{O}_j^i)= \delta_j^i = \epsilon(\bar{\mathcal{O}}_j^i),  $$
	$$S(\mathcal{O}_j^i) = \bar{\mathcal{O}}_j^i, \quad S(\bar{\mathcal{O}}_j^i) = \mathcal{O}_j^i, \quad \text{ for all } i,j \in \{1,\ldots,n\}$$
	and if any one of these conditions is met: 
	\begin{itemize}
		\item[(a)] %$H$ is commutative with $S = S^{-1}$ and there is a Hopf action $\str \colon \Ug \otimes H \to \genfd $ such that $X \str f = \langle X, f_{(1)} \rangle f_{(2)}$ and a Hopf action $\str \colon \Ug^\op \otimes H \to \genfd$ such that $Y \str f = \langle \phi(Y), f_{(2)} \rangle f_{(1)}$, where $\phi \colon \Ug^\op \to \Ug$ antiisomorphism, such that $H \sharp \Ug \cong H^\co \sharp \Ug^\op$ via $h \mapsto h, X \mapsto \sum \OO^i_j \sharp Y_j$ 
		$\genfd$ is $\mathbb{R}$ or $\mathbb{C}$, $G$ is a Lie group, $\mathcal H \subset C^\infty(G)$ as an algebra and its coalgebra structure and antipode are transpose of the group structure of $G$, and $\OO^i_j,\bar{\OO}^i_j, i,j\in \{1,\ldots,n\}$ are exactly compontents of the matrix of $\Ad$ and its inverse in the chosen basis of $\ggf$, or
		\item[(b)] $\genfd$ is $\mathbb{R}$ or $\mathbb{C}$, $G$ is a Lie group, $\mathcal H \subset C^\infty(G,e)$ as an algebra and its coalgebra structure and antipode are transpose of the group structure of $G$, and $\OO^i_j,\bar{\OO}^i_j, i,j\in \{1,\ldots,n\}$ are exactly the germs of the components of the matrix of $\Ad$ and its inverse in the chosen basis of $\ggf$, or
		\item[(c)] $\genfd$ is any field, $G$ is an affine algebraic group over $\genfd$ and $\mathcal H$ is a Hopf algebra such that $ \mathcal H \subset \OO(G)$ as a Hopf subalgebra, or
		%\item[(b)] elements $\mathcal{O}^i_j, \bar{\mathcal{O}}^i_j, i,j\in \{1,\ldots,n\}$ generate all elements of $H$ as an algebra, or
		\item[(d)] the pairing $\Ug \otimes \mathcal H \to \genfd$ is non-degenerate in the second variable, that is, $\mathcal H$ is essentially a Hopf subalgebra of the finite dual $\Ug^\circ$ of $\Ug$,
	\end{itemize}
	then $(\Ug,\btl',\lambda')$ is a braided commutative right-left Yetter--Drinfeld module algebra over $\mathcal H$, where the right action is induced by the pairing by formula $X \btl' f = \langle X_{(1)} , f \rangle X_{(2)}$ and the left coaction is defined on the Lie algebra generators by formula 
	$$\lambda'(X_j) = \sum \bar{\mathcal{O}}^i_j\otimes X_i.$$
	Likewise, $(\Ug^\op, \btl, \lambda)$ is a braided commutative right-left Yetter--Drinfeld module algebra over $\mathcal H^\co$, where the right action is induced by the pairing by formula $Y \btl f = \langle Y_{(1)} , f \rangle Y_{(2)}$ and the left coaction is defined on the Lie algebra generators by formula 
	$$\lambda(Y_j) = \sum {\mathcal{O}}^i_j\otimes Y_i.$$

\end{theorem}

\begin{proof}
	This is proven in previous theorems.
\end{proof}
This theorem covers the cases of any Hopf algebra $\mathcal H$ such that $\OO^{min}(G) \subset \mathcal H \subset C^\infty(G) \cap \mathcal R$ for a Lie group $G$ when $\genfd$ is $\mathbb{R}$ or $\mathbb{C}$, any Hopf algebra $\mathcal H$ such that $\OO^{min}(G) \subset\mathcal  H \subset \OO(G)$ when $G$ is an affine algebraic group over any field $\genfd$, and any Hopf algebra $\mathcal H$ such that $\Ug^{min}\subset\mathcal  H \subset \Ug^\circ$ for any field $\genfd$. 

\begin{corollary} Let $\mathfrak{h}$ be a finite-dimensional Leibniz algebra over any field $\genfd$. Then $U(\operatorname{Der}(\mathfrak{h}))$ is a braided commutative right-left Yetter--Drinfeld module algebra over $\OO(\Aut(\mathfrak{h}))$ and the smash product algebra $\OO(\Aut(\mathfrak{h})) \sharp U(\operatorname{Der}(\mathfrak{h}))$ is a scalar extension Hopf algebroid over base algebras $U(\operatorname{Der}(\mathfrak{h}))^\op, U(\operatorname{Der}(\mathfrak{h}))$. 
\end{corollary}

\begin{remark} 
	This is related to the example $\OO(\Aut(\mathfrak{h}))\sharp U(\mathfrak{h}_{Lie})$ in \cite{SSOAut} which embeds as an algebra into the algebra $\OO(\Aut(\mathfrak{h})) \sharp U(\operatorname{Der}(\mathfrak{h}))$ of regular differential operators on the linear affine algebraic group $\Aut(\mathfrak{h})$ of automorphisms of a Leibniz algebra~$\mathfrak{h}$ over a field $\genfd$. The proof for the example $\OO(\Aut(\mathfrak{h}))\sharp U(\mathfrak{h}_{Lie})$ uses the fact that certain functions $\GG^i_j, \bar{\GG}^i_j, i,j\in\{1,\ldots,n\}$ that correspond to a chosen basis of the Leibniz algebra $\mathfrak h$ generate the whole Hopf algebra $\OO(\Aut(\mathfrak{h}))$ and properties are all directly checked on generators without using geometrical interpretation. 
	
	Example $\OO(\Aut(\ggf)) \sharp \Ug$ from \cite{SSOAut} when $\ggf$ is Lie algebra of an affine algebraic group $G$ is connected to our examples via the adjoint representation of $G$, because the pairing of a function $f \in \OO(\Aut(\ggf))$ with an element $D \in \Ug$ in \cite{SSOAut} satisfies the formula $\langle D, f \rangle = \epsilon(D (f\circ \Ad))$. In the following diagrams all maps are Hopf algebra morphisms, with the completed version of Hopf algebra structure on $\Ug^*$ assumed, and respect the pairings of Hopf algebras with $\Ug$.
	$$
	\begin{tikzcd}
	\Ug^* & \ar{l}{} \OO(G)  &  \ar{l}{} \ar[two heads]{dl} \OO(\Aut(\ggf))  & \ar[two heads]{l} \OO(\mathrm{GL}(\ggf)) & \ar[swap]{l}{\cong} \OO(\GLnk) \\
	\Ug^{min} \ar[hook]{u}  &  \ar[two heads]{l} \OO^{min}(G) \ar[hook]{u}  & \text{ for }\ggf = \gL & &
	\end{tikzcd}
	$$
	$$
	\begin{tikzcd}
	(\Ug^*)^\co & \ar{l}{} \OO(G)  &  \ar{l}{} \ar[two heads]{dl} \OO(\Aut(\ggf))  & \ar[two heads]{l} \OO(\mathrm{GL}(\ggf)) & \ar[swap]{l}{\cong} \OO(\GLnk) \\
	(\Ug^{min})^\co \ar[hook]{u}  &  \ar[two heads]{l} \OO^{min}(G) \ar[hook]{u}  & \text{ for }\ggf = \gR & &
	\end{tikzcd}
	$$
	
\end{remark}
\section{Hopf algebroids} We present a Hopf algebroid structure, that is a scalar extension of $\UgL$ with $\OO(G)$, over base algebras $\UgR, \UgL$, on algebra 
$$\operatorname{Diff}(G) \cong  \OO(G) \sharp \UgL \cong \OO(G)^\co \sharp \UgR \cong \UgR \sharp \OO(G) \cong \UgL \sharp \OO(G)^\co,$$  and additionally another Hopf algebroid structure, a scalar extension of $\UgR$ with $\OO(G)^\co$, over base algebras $\UgL, \UgR$ on the same algebra. Similarly, for any Hopf algebra such that $\OO^{min}(G) \subset \mathcal H \subset \OO(G)$ for affine algebraic group $G$ over a field $\genfd$, or any Hopf algebra such that $\OO^{min}(G) \subset \mathcal H \subset C^\infty(G) \cap \mathcal R$ for a Lie group $G$ over $\mathbb R$ or $\mathbb C$, Hopf algebroid structures related to the ones above exist on algebra 
$$\KK \cong \UgR \sharp\mathcal  H \cong\mathcal  H \sharp \UgL \cong\mathcal  H^\co \sharp \UgR \cong  \UgL \sharp \mathcal H^\co \subset \operatorname{Diff}(G).$$

We also present related structures of a Hopf algebroid on $\Ug^\circ \sharp \Ug$  over base algebras $\Ug^\op, \Ug$ and, similarly, on $\Ug^\circ \sharp \Ug$ over base algebras $\Ug,\Ug^\op$.   Similarly, for any Hopf algebra such that $\Ug^{min} \subset \mathcal H \subset \Ug^\circ$, Hopf algebroid structures related to the ones above exist on the algebra 
$$\mathcal  H \sharp \Ug \cong\Ug^\op\sharp\mathcal  H^\co \subset \Ug^\circ \sharp \Ug .$$ 
When $G$ is a Lie group over $\mathbb R$ or $\mathbb C$ and $\ggf$ is its Lie algebra, the smash product algebra $\Ug^* \sharp \Ug$ is isomorphic to the algebra $ \operatorname{Diff}^\omega(G,e)\sharp \Ug \cong J^\infty(G,e)\sharp \Ug$ of formal differential operators around the unit of  $G$.

\subsection{Diagram of scalar extensions}
The following diagram shows the relationship of different smash product algebras containing $\Ug$ as a factor and induced by the pairing of $\Ug$ with Hopf algebras. The general formulas for scalar extension are specialized for generators of these algebras in Subsection~\ref{subs:lasttable}. For the left invariant case, we have
$$
\begin{tikzcd}
\OO(\Aut(\ggf)) \sharp U(\ggf) \ar[two heads]{d}   &  &  
\\
\OO^{min}(G) \sharp \Ug \ar[two heads]{d}  \ar[hook]{r} &  \OO(G) \sharp U(\ggf) \ar{d} &\text{ for } \ggf = \gL
\\
\Ug^{min} \sharp U(\ggf)  \ar[hook]{r} & \Ug^{\circ} \sharp U(\ggf)\ar[hook]{r} & \Ug^*\sharp \Ug,
\end{tikzcd}
$$
%We have a structure of a Hopf algebroid over $\Ug^\op, \Ug$ on the following algebras:
%\begin{itemize}
%	\item[(i)] Minimal scalar extension $\Ug^{min} \sharp \Ug  \subseteq \Ug^{*} \sharp U(\ggf)$, for a Lie algebra $\ggf$
%	%	\item[(ii)] Heisenberg double with the reduced dual $\Ug^{\circ} \sharp U(\ggf)$, for a Lie algebra $\ggf$
%	\item[(ii)] Minimal algebra $\OO^{min}(G) \sharp U(\ggf) \subseteq \mathrm{Diff}(G)$, for a Lie group $G$ 
%	\item[(iii)] Algebra $\OO(G) \sharp U(\ggf)$ of regular differential operators on affine algebraic group~$G$
%	\item[(iv)] Finite dual Heisenberg double $\Ug^\circ\sharp \Ug$ for a Lie algebra $\ggf$
%\end{itemize}
and for the right invariant case, we have
$$
\begin{tikzcd}
\OO(\Aut(\ggf))^\co \sharp U(\ggf) \ar[two heads]{d} & & 
\\
\OO^{min}(G)^\co \sharp \Ug \ar[two heads]{d}  \ar[hook]{r} &  \OO(G)^\co \sharp U(\ggf) \ar{d}  & \text{ for } \ggf = \gR 
\\
\Ug^{min} \sharp U(\ggf)  \ar[hook]{r} &  \Ug^{\circ} \sharp U(\ggf) \ar[hook]{r} &  \Ug^*\sharp \Ug.
\end{tikzcd}
$$

\subsection{Formulas in generators and structure constants}
\label{subs:lasttable}
%Here we give formulas for all scalar extensions appearing in this article. Each smash product algebra $H \sharp \UgL$ can be presented as an algebra in four different ways, as
%$$
%H \sharp \UgL \cong \UgR \sharp H \cong H^{\mathrm{co}} \sharp \UgR \cong \UgL \sharp H^{\mathrm{co}},
%$$
%and we give the structure map formulas for each such presentation. In the table we give explicit formulas for $H= \UgL^{min}, H^\co \cong \UgR^{min}$, as these are images of both algebras $\OO^{min}(G)$ and $\OO(\Aut(\ggf))$ under mappings to duals $\UgL^*, \UgR^*$ induced by pairings with $\UgL, \UgR$ respectively; to obtain the formulas for $H = \OO^{min}(G)$ and $H = \OO(\Aut(\ggf))$, the generators that are the matrix components of $\mathcal{U}$ are to be replaced with the correspoding generators of the matrix $\mathcal{O}$ or $\mathcal{G}$, respectively.
%
%Induced by the pairing, we have mappings $H \to \UgL^*$, $H^{\mathrm{co}} \to \UgR^*$.  
%Denote the comultiplication on  $H$ in Sweedler notation by $$\Delta(f) = f_{(1)} \otimes f_{(2)}.$$ Therefore, the comultiplication on  $H^{\mathrm{co}}$ is here $\Delta^{\mathrm{co}}(f) = f_{(2)} \otimes f_{(1)}$. 
%

\subsubsection{Scalar extension of $\UgR$ with $\mathcal H^\co$} Denote $R := \UgR$ and $L := \UgL$. Let $\phi \colon \UgL \to \UgR$ be the antiisomorphism defined on generators by $X_i \mapsto Y_i$ for $ i\in\{1,\ldots,n \}$. We here present formulas for the structure of a Hopf algebroid over the base algebras $\UgL, \UgR$ on the smash product algebras 
%$$ \UgL \sharp \UgR^* \cong \UgR^* \sharp \UgR.$$ The formulas are the same,  as is explained above, to the formulas for the smash product algebras 
$$ \UgL \sharp \mathcal H^{\mathrm{co}} \cong \mathcal H^{\mathrm{co}} \sharp \UgR$$
for any algebra $\mathcal H$ such that: $\UgL^{min} \subset \mathcal H \subset \UgL^\circ$ or $\OO^{min}(G) \subset \mathcal H \subset \mathcal \OO(G) $.

Since $(\UgR, \btl, \lambda)$ is a right-left Yetter--Drinfeld module algebra over $\mathcal H^{\co}$, by the results in \cite{stojic}, $\mathcal H^\co \sharp \UgR \cong \UgL \sharp \mathcal H^\co$ is a scalar extension Hopf algebroid over base algebras $\UgL, \UgR$, with the following structure map formulas.

\begin{align*}
&\text{In terms of } 	\UgL \sharp \mathcal H^\co:& &  \text{In terms of }                \mathcal   H^\co \sharp \UgR:\\
&&& \\
& \textit{left bialgebroid over } \UgL & 
& \text{left bialgebroid over }  \UgL \\ 
& \alpha_{L}(X_j) = X_j& 
& \alpha_{L}(X_j) = X_j = \textstyle\sum_i {\UU}_j^i \sharp Y_i   \\
& \beta_{L}(X_j) = \textstyle\sum_i \bar\UU_j^i \cdot X_i - \textstyle\sum_i  C_{ij}^i     	 & 
& \beta_{L}(X_j) =    Y_j - \textstyle\sum_i C_{ij}^i \\
& \Delta_{L}(X \sharp f) = X \sharp f_{(2)} \otimes_{L} 1\sharp f_{(1)}    &
& \Delta_{L} (f \sharp Y) =  f_{(2)} \sharp 1 \otimes_{L} f_{(1)} \sharp Y  \\
& \epsilon_{L}(X \sharp f) =\epsilon(f) X &   
&  \epsilon_{L}(f\sharp Y_j)  = f \btr X_j - \epsilon(f) \textstyle\sum_i C_{ij}^i \\ 
%	\end{align*} 
%	
%	\begin{align*}
%	&	\UgL \sharp H^\co& &                    H^\co \sharp \UgR\\
&&& \\
& \text{right bialgebroid over } \UgR &
&  \textit{right bialgebroid over }  \UgR \\   
& \alpha_R(Y_j)= Y_j = \textstyle\sum_i \bar\UU_j^i \cdot X_i  &  
&                   \alpha_R(Y_j) =  Y_j   \\
& \beta_R(Y_j) = X_j     &   
&  \beta_R(Y_j) = X_j = \textstyle\sum_i {\UU}_j^i \sharp Y_i  \\
& \Delta_R(X \sharp f) = X\sharp f_{(2)} \otimes_{R} 1 \sharp f_{(1)}      &
&  \Delta_R(f \sharp Y) = f_{(2)} \sharp 1 \otimes_{R} f_{(1)} \sharp Y   \\
& \epsilon_R(X_j \sharp f) = Y_j \btl f  &
&  \epsilon_R(f \sharp Y) = \epsilon(f)Y   
%\hline
%\end{align}
%
%\begin{align*}{ll}
\\
\\
& \text{antipode}  & & \text{antipode}  \\ 
& \tau (X_j \sharp f) =  S(f) \cdot Y_j  &
&  \tau (f \sharp Y_j) = ( X_j  + \textstyle\sum_i C_{ij}^i) \cdot S(f) \\
& \tau^{-1} (X_j \sharp f) = S^{-1}(f)(Y_j - \textstyle\sum_i C_{ij}^i)  & 
& \tau^{-1} (f \sharp Y_j) = X_j \cdot S^{-1}(f)\\
\\
& \text{Here }\phi = \epsilon_R \circ \alpha_L,\ \phi^{-1} = \epsilon_L \circ \beta_R,\ & & \phi \colon \UgL \to \UgR \text{ antiisomorphism.} 
\end{align*}
In the formulas for the antipode and its inverse, and for the counit, on the left side we have written $Y_j$ as a shorthand for its smash product algebra form on the left $$Y_j = \sum_i\bar{\mathcal{U}}^i_j \cdot X_i \in \UgL\sharp\mathcal H^\co,$$ and on the right side we have written $X_j$ as a shorthand for its smash product  algebra form on the right,  $$X_j = \sum_i \mathcal{U}^i_j\sharp  Y_i \in \mathcal H^\co \sharp \UgR.$$ The actions appearing in formulas for the counits are the ones that define the corresponding smash product algebras; they are: the right Hopf action $$\btl \colon \UgR \otimes \mathcal H^\co \to \UgR$$ defined previously from the pairing, and the left Hopf action $$\btr \colon \mathcal H^\co \otimes \UgL \to \UgL$$ defined as $f  \btr X = \phi^{-1}(\phi(X)\btl S(f))$, for $f\in \mathcal H^\co$ and $X \in \UgL$.

The previous formulas are then translated to the other two smash product algebra forms of the algebra $\KK$, through the algebra isomorphisms %$(\KK, \alpha_L, \beta_L, \Delta_L, \epsilon_L, \alpha_R,\beta_R,\Delta_R,\epsilon_R, \tau )$,
$$\KK \cong\mathcal H \sharp \UgL \cong \UgR \sharp\mathcal H \cong \mathcal H^{\mathrm{co}} \sharp \UgR \cong \UgL \sharp\mathcal H^{\mathrm{co}}.
$$
In the next table there appear the same two actions as in the table above; they are not the actions that provide these two new smash product algebra structures.
\begin{align*}
&\text{In terms of } \mathcal	H\sharp \UgL: &
&  \text{In terms of }                 \UgR \sharp\mathcal H: \\
&&& \\
& \textit{left bialgebroid over } \UgL &
& \text{left bialgebroid over }  \UgL\\ 
& \alpha_{L}(X_j) = X_j&
& \alpha_{L}(X_j) = X_j = \textstyle\sum_i {\UU}_j^i \cdot Y_i   \\
& \beta_{L}(X_j) =  \textstyle\sum_i \bar\UU_j^i\sharp X_i - \textstyle\sum_i C_{ij}^i         &
& \beta_{L}(X_j) =    Y_j - \textstyle\sum_i C_{ij}^i \\
& \Delta_{L}(f \sharp X) = f_{(2)} \sharp 1 \otimes_{L}  f_{(1)} \sharp X    &
&   \Delta_{L} (Y \sharp f) = Y \sharp f_{(2)}  \otimes_{L} 1 \sharp f_{(1)}  \\
& \epsilon_{L}(f \sharp X) = f \btr X & 
&   \epsilon_{L}(Y_j\sharp f)  = \epsilon(f)(X_j - \textstyle\sum_i C^i_{ij})\\ 
%	\end{align*} 
%	
%	\begin{align*}
%	&	\UgL \sharp H^\co& &                    H^\co \sharp \UgR\\
&&& \\
& \text{right bialgebroid over } \UgR &
&  \textit{right bialgebroid over }  \UgR \\   
& \alpha_R(Y_j)= Y_j = \textstyle\sum_i \bar\UU_j^i \sharp X_i  &
&  \alpha_R(Y_j) =  Y_j   \\
& \beta_R(Y_j) = X_j     & 
&  \beta_R(Y_j) = X_j = \textstyle\sum_i {\UU}_j^i \cdot Y_i  \\
& \Delta_R(f \sharp X) = f_{(2)} \sharp 1 \otimes_{R} f_{(1)} \sharp X      &
&  \Delta_R(Y \sharp f) = Y \sharp f_{(2)} \otimes_{R} 1 \sharp f_{(1)}    \\
& \epsilon_R(f \sharp X_j) = \epsilon(f) Y_j  &
&  \epsilon_R(Y \sharp f) = Y \btl f   
%\hline
%\end{align}
%
%\begin{align*}{ll}
\\
\\
& \text{antipode}  &
& \text{antipode}  \\ 
& \tau (f \sharp X_j) =  Y_j \cdot S(f) &
&  \tau (Y_j \sharp f) = S(f) (X_j + \textstyle\sum_i C_{ij}^i) \\
& \tau^{-1} (f \sharp X_j) = (Y_j - \textstyle\sum_i C_{ij}^i)  S^{-1}(f) &
& \tau^{-1} (Y_j \sharp f) = S^{-1}(f) \cdot X_j  \\
\\
& \text{Here }\phi = \epsilon_R \circ \alpha_L,\ \phi^{-1} = \epsilon_L \circ \beta_R,\ &
& \phi \colon L \to R \text{ antiisomorphism.} 
\end{align*}

\subsubsection{Scalar extension of $\UgL$ with $\mathcal H$} On the other hand, let now $L' := \UgR$ and $R' := \UgL$. Let $\psi \colon \UgR \to \UgL$ be the antiisomorphism defined on generators by $Y_i \mapsto X_i$. It is the inverse of the former antiisomorphism $\phi$. Since $(\UgL, \btl', \lambda')$ is a right-left Yetter--Drinfeld module algebra over $\mathcal H$, by the results in \cite{stojic}, $\mathcal H \sharp \UgR \cong \UgL \sharp \mathcal H$ is a scalar extension Hopf algebroid over base algebras $\UgR, \UgL$. We have chosen to write formulas for the previous case, since here the left bialgebroid is over $\UgR$ and the right bialgebroid is over $\UgL$. Formulas for this case are analogous, see the following subsection.
\subsubsection{Mirror picture}
One could also view the algebra $ \Diff(G)$ as if acting from the right on functions, and obtain mirror smash product algebras to the ones here. That would produce a scalar extension of $\UgR$ with $\mathcal H$ over base algebras $\UgL, \UgR$ and a scalar extension of $\UgL$ with $ \mathcal H^\co$ over base algebras $\UgR,\UgL$. Such mirror construction corresponds to the one used to derive a version of a completed Hopf algebroid structure on a Lie algebra type noncommutative phase space in \cite{halg}. 
%\section{Hopf algebroid formulas}\label{sec:HA}

\begin{remark} The smash product algebra $\Ug^{*} \sharp \Ug$, which can be identified with the algebra of formal differential operators  $\mathrm{Diff}^\omega(G,e)\cong J^\infty(G,e)\sharp \Ug$, is defined as an internal Heisenberg double in the category $(\ipV,\tilde\otimes,\genfd)$ of filtered cofiltered vector spaces, introduced in the dissertation of one of the authors \cite{stojicphd}. It is there shown to be an internal Hopf algebroid \cite{stojicint} in the same category.  Hopf algebroids $\Ug^{min} \sharp \Ug$ and $\Ug^\circ\sharp \Ug$ defined here embed in this internal Hopf algebroid.
\end{remark}
%and the following internal Hopf algebroid over $\Ug_\mathbb{C}$
%\begin{itemize}
%	\item[(vi)] $\Rep(G) \sharp \Ug_\mathbb{C}$, for a connected compact Lie group $G$.
%\end{itemize}

\subsection{Formulas in short}

\subsubsection{Scalar extension of $\UgR$ with $\mathcal H^\co$}
In the scalar extension of $\UgR$ with $\mathcal H^\co$, the algebra $\KK \cong \mathcal H^\co \sharp \UgR \cong \UgL \sharp \mathcal H^\co$ is a left bialgebroid over $\UgL$ and a right bialgebroid over $\UgR$ with the source map and the target map satisfying:
$$\alpha_L(X) = X, \quad \alpha_R(Y)=Y, \quad  \beta_L(X_j)= Y_j - \sum_i C_{ij}^i,   \quad \beta_R(Y_j) = X_j,$$
the comultiplication and the counit satisfying:	
$$\Delta_L(X) = X\otimes_L 1, \quad  \Delta_L(Y) = 1\otimes_L Y, \quad \Delta_R(X) = X\otimes_R 1, \quad \Delta_R(Y) = 1\otimes_R Y,$$ 
$$\Delta_L(f)=  \Delta_R(f)= \Delta^{\co}(f),$$
$$\epsilon_L(X \cdot f) = \epsilon(f) X, %\quad \epsilon_L(f \cdot X) = f\btr X, 
\quad \epsilon_R(f \cdot Y) = \epsilon(f)Y, %, \quad \epsilon_R(Y \cdot f ) = Y \btl f
$$
comprising a Hopf algebroid with the antipode satisfying: 
$$\tau(X_j) = Y_j, \quad \tau(f) = S(f).$$ 
Here $X$ and $Y$ denote generic elements of $\UgL$ and $\UgR$ respectively as subsets of $\mathcal K$, and $\{X_1,\ldots, X_n\}$ and $\{Y_1,\ldots,Y_n\}$ are the chosen bases for $\gL$ and $\gR$, respectively, that coincide at the unit of $G$. These structure maps are not dependent on the choice of the basis. 

The actions $\btr \colon\mathcal  H^\co \otimes \UgL \to \UgL$ and $\btl \colon \UgR \otimes\mathcal  H^\co \to \UgR$ that define the smash product algebras $\UgL \sharp\mathcal  H^\co$ and $\mathcal H^\co \sharp \UgR$ respectively satisfy
$$f \btr X = \epsilon_L(f \cdot X), \quad Y \btl f = \epsilon_R(Y \cdot f ).$$
\subsubsection{Scalar extension of $\UgL$ with $\mathcal H$}
In the scalar extension of $\UgL$ with $\mathcal H$, the same algebra $\KK \cong \mathcal H \sharp \UgL \cong \UgR \sharp \mathcal H$ has the structure of a left bialgebroid over $\UgR$ and a right bialgebroid over $\UgL$ with the source map and the target map satisfying:
$$\alpha'_{L}(Y) = Y, \quad \alpha'_{R}(X)=X, \quad \beta'_{L}(Y_j)= X_j + \sum_i C_{ij}^i, \quad \beta'_{R}(X_j) = Y_j,$$
the comultiplication and the counit satisfying:	
$$ \Delta'_{L}(Y) = Y\otimes_{L'} 1, \quad \Delta'_{L}(X) = 1\otimes_{L'} X, \quad \Delta'_{R}(Y) = Y\otimes_{R'} 1, \quad \Delta'_{R}(X) = 1\otimes_{R'} X $$ 
$$\Delta'_{L}(f)=  \Delta'_{R}(f)= \Delta(f)$$	
$$\epsilon'_{L}(Y \cdot f) =  \epsilon(f) Y, \quad \epsilon'_{R}(f \cdot X) = \epsilon(f)X$$
comprising a Hopf algebroid with the antipode map satisfying  $$\tau'(Y_j) = X_j, \quad \tau'(f) = S(f).$$
Since here we have $S(f) = S^{-1}(f)$, it follows that $\tau$ is inverse to the above $\tau'$. We have that $\alpha_L' = \alpha_R$, $\alpha_L = \alpha_R'$, that $\epsilon_R \circ \beta_L$ is inverse to $\epsilon_R' \circ \beta_L'$, and that $\epsilon_L \circ \beta_R$ is inverse to $\epsilon_L' \circ \beta_R'$. 

\section*{Relation to the previous work} Examples $\OO^{min}(G)$, $\Ug^{min}$ and $\Ug^\circ$ are first presented as examples as part of Chapter 9 of the first author's dissertation  \cite{stojicphd} of title \emph{Completed Hopf algebroids}. Matrix $\UU$ there arose naturally while considering a Hopf algebroid structure on the  Heisenberg double $\Ug^* \sharp\Ug$ internally in the category of filtered cofiltered vector spaces introduced in the dissertation. Components of $\UU$ were precisely coefficients appearing in a certain formal sum that is an infinite version of Lu's formula  \cite{Lu} for coaction. Matrix $\OO$ was introduced geometrically there to satisfy the same properties as the corresponding matrix in the example $\OO(\Aut(\ggf))$ \cite{SSOAut} of the second author.  This material is then expanded here to include $\OO(G)$ over any field $\genfd$ and $\Ug^\circ$ without the use of the results that need completions.

\end{document}